\UseRawInputEncoding
\documentclass[10pt]{amsart}
\usepackage{graphicx,amssymb,amsmath}
\newtheorem{thm}{Theorem}[section]
\newtheorem{cor}[thm]{Corollary}
\newtheorem{lem}[thm]{Lemma}
\newtheorem{prop}[thm]{Proposition}

\theoremstyle{definition}
\newtheorem{defn}[thm]{Definition}

\theoremstyle{remark}
\newtheorem{rem}[thm]{Remark}
\newtheorem*{notation}{Notation}

\numberwithin{equation}{section}

\def\Z   {{\bf Z}}

\def\R   {{\bf R}}
\def\C   {{\bf C}}
\def\P   {{\bf P}}

\def\calK      {{\mathcal K}}

\def\O           {{\mathcal O}}
\def\calE       {{\mathcal E}}

\def\Kahler    {K\"ahler\ }
\def\CY        {Calabi--Yau\ }

\def\vol         {{\rm vol}}

\def\Z   {{\bf Z}}

\def\R   {{\bf R}}
\def\C   {{\bf C}}
\def\P  {{\bf P}}

\newcommand{\beas}{\begin{eqnarray*}} 
\newcommand{\eeas}{\end{eqnarray*}} 

\begin{document}

\title{The Topology of Calabi--Yau threefolds with Picard number three}


%    author one information

\author{P.M.H. Wilson} 
\address{Department of Pure Mathematics, University of Cambridge,
16 Wilberforce Road, Cambridge CB3 0WB, UK}
\email {pmhw@dpmms.cam.ac.uk}

\date{15 May 2023}

%%%%%%%%%%%%%%%%%%%%%%%%%%%%%%%%%%%%%%%%%%%%%%%%%%%%%%%%%%%%%%%%
\begin{abstract}
We ask about  the simply connected compact smooth  6-manifolds which can support structures of Calabi--Yau threefolds.  In particular, we study the interesting case of Calabi--Yau threefolds $X$ 
with the second betti number $b_2 (X)  =3$.  We have  a 
cubic form $D \mapsto D^3$ on $H^2 (X, {\bf Z})$ given by cup-product, a linear form 
 $D \mapsto D\cdot c_2(X)$ given by the second Chern class, and the integral middle cohomology $H^3 (X, \Z)$, and if 
 $X$ is simply connected with torsion free homology  this information determines 
precisely  the 
diffeomorphism class of the underlying 6-manifold by a result of Wall.  For simplicity, we assume that the cubic form defines a smooth real elliptic curve whose Hessian is 
irreducible.  Under a further 
relatively mild assumption that there are no non-movable surfaces $E$ on $X$ with $1\le E^3 \le 8$,
 we prove that the real elliptic curve must have two connected components rather 
than one, and that the \Kahler cone  is contained in the open positive cone  on the bounded component;  we show moreover that the linear form $c_2$ is also non-negative on this cone.  Using Wall's result, for any given third betti number
we therefore have an abundance of examples of smooth compact oriented 6-manifolds which support no Calabi--Yau structures, both in the cases when the cubic defines a real elliptic curve with one or two connected components. Moreover, except possibly if $c_2$ vanishes at a real inflexion point of the elliptic curve,
even when  \CY structures do occur under the above conditions, there will be only a bounded family of them which are not birationally elliptic.

 \end{abstract}

%%%%%%%%%%%%%%%%%%%%%%%%%%%%%%%%%%%%%%%%%%%%%%%%%%%%%%%%%%%%%%%%

\maketitle

\section*{Introduction}
\rm  About 35 years ago, the author was asked by two colleagues in physics what restrictions there were on the cup-product cubic form on $H^2 (X, \Z)$ given by $D \mapsto D^3$ for a Calabi--Yau threefold $X$.  Until now, we were only able to make rather weak statements in response to this question, but recent work now enables us to give an answer 
with more content.  The question becomes interesting when $b_2 (X) \ge 3$, and in this paper we concentrate on the case $b_2(X) =3$; the intervening three decades suggested the following Motivating Question to the author:
\vspace{2 mm}

\noindent \bf Motivating Question. \rm Is the following true?  If $X$ is a  
\CY threefold with $b_2 (X)=3$ and the cup-product cubic form is smooth, therefore defining an elliptic curve with real locus 
$C \subset \P^2 (\R )$, then either $C$ has two connected components and the \Kahler cone $\calK (X)$ is contained in the cone on the bounded component of $C$ on which the cubic is positive, or $C$ has one real  component and the Hessian curve of the cubic 
is reducible, consisting of three non-concurrent lines.

\begin{rem}  The cases with reducible Hessian do commonly  occur: for one possibility we consider the desingularization of 
a quintic in $\P ^4 (\C )$ with two singularities, each of which is analytically a cone on a del Pezzo surface, and for another we take the desingularization of a hypersurface of bidegree $(3,3)$ 
in $\P^2 \times \P^2$ with one such singularity.  The Hessian curve in the first case consists of three real lines and in the second case one real and two complex lines.  In the real coordinates introduced in Section 1, the corresponding invariant $k$ is $0$ and $-2$ respectively.  

The answer however to the question as posed is in fact no, as there are some known examples 
%(all admitting contractions of  surfaces on $X$ to a curve $\P^1$ of singularities) 
where the real elliptic curve $C$ has one component and the Hessian is smooth.  Three such examples arise from hypersurfaces in weighted projective space and are included in the 
list given in Appendix C of \cite{HLY}.
 If we consider for instance a quasi-smooth hypersurface $X_{10}$ of degree 10 in weighted projective space $\P(1,1,2,3,3)$, 
 there is a curve (isomorphic to $\P^1$) of singularities.  Taking a crepant resolution, we obtain a \CY threefold $X$ with $b_2 (X) =3$.
 %There are however only 12 
 %non-simply connected toric hypersurface \CY threefolds of which three have $b_2 =3$, which are listed in Section 2 of \cite{BK}; the given threefold is not one of these and hence is simply connected.
 With respect to a suitable 
 basis of $H^2(X, \Z )$ the cubic form is given explicitly in Appendix C of \cite{HLY} as
 $$15x^3 +60x^2y +30 x^2z +78 xy^2 +78 xyz + 18xz^2 +32y^3 +48 y^2z +18yz^2.$$
   This may be checked to be an example where the real elliptic curve has one connected component, and the Hessian is also smooth.  
   \end{rem}
We recall that the \it movable cone \rm of a smooth threefold $X$ is the closure of the cone in $H^2(X, \R)$ generated by the classes of mobile divisors, and it follows from 
 Lemma 3.2 of \cite{WilBd} that when $b_2(X)=3$,   
the Hessian is non-negative on the movable cone.
As we saw in the earlier paper \cite{WilBd}, the \it rigid non-movable \rm surfaces $E$ on 
a \CY threefold $X$, which were defined as the irreducible surfaces on $X$ that deform with any small deformation of the complex structure on $X$ but for which no multiple moves in the threefold (see Section 2 of \cite{WilBd} for a discussion of these) play a crucial role in understanding possible \CY structures on a compact 6-manifold.  
If we impose the extra condition that there are no rigid non-movable surfaces on $X$ with $E^3 >0$, then the answer to the above question is yes.

Recall the result of Wall (\cite{Wall}, Theorem 5) that under an assumption ($H$) that the compact 
simply connected oriented 6-manifolds studied have torsion free homology and class $w_2 (M) =0$ (the latter assumption holding if $M$ supports a \CY structure), then the diffeomorphism classes of compact oriented manifolds $M$ satisfying ($H$) correspond bijectively to isomorphism classes of invariants:

two free abelian groups $H$ and $G$ (corresponding to $H^i (M, \Z)$ for $i=1,2$) with the rank of $G$ being even,

a symmetric trilinear map $\mu : H \times H \times H \to \Z$,

a homomorphism $p_1 : H \to \Z$,

subject to: for all $x, y \in H$,
$$\mu (x,x,y) \equiv \mu (x,y,y)\quad (\hbox{\rm mod\ } 2)\quad \hbox{\rm and} \quad p_1 (x) \equiv 4\mu (x,x,x) \quad (\hbox{\rm mod\ } 24).$$
We shall be interested in the case when there is a \CY structure $X$ on the manifold, in which case the linear form on the given free abelian group $H$ may be identified as $p_1 (X) = -2c_2(X)$.  Wall first constructs the relevant 6-manifold $M_0$ when $G=0$, and he then forms a connected sum of $M_0$ with $b_3/2$ copies of $S^3 \times S^3$, and so at the smooth level the information on $H^3 (M, \Z)$ and the invariants on $H^2 (M, \Z )$ are independent.  
This is  not true for \CY threefolds, since 
at least for $b_2 (X)\le 2$ we know that  the invariants on $H^2(X, \Z)$ 
(non-effectively) bound $b_3(X)$ (for $b_2 =1$ this follows using Hilbert schemes, and see \cite{WilBd} for the $b_2 =2$ case).
When $M$ supports an almost complex structure with $c_1 =0$, it is unique up to homotopy by Theorem 9 of \cite{Wall}.

If $E$ is a rigid non-movable surface with $E^3>0$ on a Calabi--Yau threefold $X$, then by the results of Section 2 from \cite{WilBd}, any such $E$ would have $E^3 \le 9$, and if $E^3 =9$  there would be a contraction of $E \cong \P^2$ to a point and the Hessian would be reducible.  
Moreover by 
applying (maybe a number of times) formula ($1'$) in Section 2 
of \cite{WilBd}, 
we deduce under the further assumption that the Hessian curve is non-singular
 that \it  either \rm $E^3 =1$ and $c_2(X)\cdot E = -2$, \it or \rm $1 \le E^3 \le 8$ and $c_2(X)\cdot E = 12 - 2E^3$.  
 Only when the integral cubic and linear forms $\mu$ and $c_2 = -p_1/2$  represent one of the above nine pairs of values at some point of $\Z^3$ with the correct index could there exist a rigid non-movable surface $E$ with $E^3 >0$.
 If $\mu$ and $p_1$ satisfy the above congruences (for instance arising from some \CY threefold), then so too do integral multiples of them, and taking an appropriate integral multiple will rule out the possibility of classes with the above pairs of invariants --- for instance scaling both $\mu$ and $p_1$ by a factor of 4 will suffice, since then $E^3$ could only be 4 or 8, but $c_2(X)\cdot E$ would be a multiple of 8.

\begin{thm}  Let $X$ be a \CY threefold with $b_2 (X)=3$, where the cubic form $F$ defines a smooth cubic curve with irreducible Hessian and there are no rigid non-movable surfaces $E$ on $X$ with $E^3>0$.  Then the real elliptic curve $C \subset \P^2 (\R )$ determined by $F$ has two connected components and the \Kahler cone of $X$ is contained in the cone on the interior of the bounded component on which the cubic is positive.  Moreover the linear form $c_2$ is non-negative on the positive  cone 
in $H^2(X, \R)$ on the bounded component.
\end{thm}

In particular, we note using Wall's result that, for any given even $b_3 >0$, 
we  have an abundance of examples of smooth compact oriented 6-manifolds which support no Calabi--Yau structures, both in the case when the cubic defines a real elliptic curve with one component and in the case of two components --- for the latter we may also need to choose the linear form appropriately, so that it is negative somewhere on the 
positive cone on the bounded component.  

It is an open question for $b_2 (X) >3$, assuming the cubic form is smooth with irreducible Hessian and $X$ contains no rigid non-movable surfaces $E$ with $E^3 >0$, 
whether the corresponding  real cubic hypersurface must contain two components with the \Kahler cone of $X$  contained in the cone on the interior of the bounded component on which the cubic is positive --- by Remark 4.7 of \cite{WilBd}, it is true when there are no rigid non-movable surfaces. 

% The proof we give below in the case $b_2 (X) =3$ depends crucially in a number of places on classical facts concerning elliptic curves.

 The last sentence of the above theorem has consequences for boundedness questions, which we explore in the final Section of the paper.   If $c_2$ were numerically trivial 
then, by a well-known result due to S.-T. Yau \cite{Yau},  $X$ is an \'etale quotient of an abelian threefold and so is not  simply connected; thus $c_2 =0$ defines a projective line
in $\P^2(\R )$.
 
\begin{thm}  Given a smooth integral ternary cubic $F$ with irreducible Hessian and an integral linear form $c_2$ not vanishing at any of the three real inflexion points of the elliptic 
curve, the \CY threefolds $X$ with these invariants on $H^2(X, \Z )$
 which do not contain rigid non-movable surfaces $E$ with $E^3 >0$ form a bounded (maybe empty) 
family, except possibly for those which are birationally elliptic.  In particular the family of such threefolds is birationally bounded.
\end{thm}

\section{Components of the  positive index cone for real ternary cubics}

In \cite{WilBd}, a central role was played by the \it positive index cone \rm corresponding to the cubic on $H^2 (X, \Z ) = \Z ^\rho$, 
namely the real classes $L$ for which $L^3 >0$ and the quadratic form given by $D \mapsto L\cdot D^2$ has index $(1,\rho -1)$.
It was shown in \cite{WilBd} that this cone has finitely many connected components, which will be convex and one of which willl contain the \Kahler cone.  If the real cubic hypersurface is smooth with two components, 
then there is a component of the positive index cone which is the positive cone on the bounded component; otherwise the Hessian will vanish at some points of the boundary of any other component of the positive index cone.  In the absence of rigid non-movable surfaces on $X$, it was remarked in Remark 4.7 of \cite{WilBd} (noting a misprint where `small locus' should read `smooth locus') that if a component of the positive index cone contains the \Kahler cone, then its boundary cannot contain points at which the cubic is strictly positive; we observed in \cite{WilBd} that under the assumption that $X$ is general in moduli, any such point would represent a big movable real divisor, which  by the theory of \it volume \rm of divisors would also be true for all nearby points.  This would produce a contradiction since the Hessian is strictly negative at some nearby points but is always 
non-negative on movable classes.  In the case $b_2 (X) =3$ we have a precise description of the components of the positive index cone, which enables us to extend this idea to 
the case when there are rigid non-movable surfaces on $X$, yielding a proof of Theorem 0.2.  

\begin{defn}  A component of the positive index cone is called \it hybrid  \rm if there are points on the boundary at which the cubic vanishes and the Hessian is strictly positive and 
points on the boundary at which the Hessian vanishes and the cubic is strictly positive.
\end{defn}

It turns out that the crucial case to rule out is the \Kahler cone being contained in a  hybrid  component of the positive index cone; 
in summary, this will be achieved when $b_2(X)=3$ as follows:

\smallskip
(1) Using the classical theory of elliptic curves, if $E_1, E_2, \ldots$ denote the rigid non-movable surface classes in $H^2 (X, \Z )$  (in particular  having Hessian non-negative there), then if $E_i \le 0$ for all $i$, the arguments in Sections 2, 3 and 4 of this paper show that part of the boundary of any hybrid component consists of points $D$ at which the cubic is positive and the Hessian vanishes, but 
which are not visible in the sense of convex geometry (with respect to the given hybrid component) from any of the $E_i$, and for which in addition $E_i \cdot D^2 >0$ for all $i$.  
In fact this part will be proved under the weaker assumption described below that all the $E_i$ lie in the extended upper half-space.
Although appearing rather technical, the  mathematical content of the proof given is essentially just   the classical theory of the Steinerian involution on the Hessian curve (since that tells us about the points $D$ on the boundary of the given hybrid component at which the Hessian vanishes and for which $E_i \cdot D^2$ is positive  for all $i$).  This proof 
does not depend on further properties of the threefold $X$ with $b_2  (X) =3$, such as being Calabi--Yau.

(2) If such a hybrid component  were to contain the \Kahler cone, under the assumption that 
the \CY threefold $X$ is general in moduli, it follows from Proposition 4.1 of \cite{WilBd} that a general such boundary class $D$ would represent a real pseudo-effective divisor
(the limit of rational big divisors); 
this follows essentially by Riemann--Roch and Kodaira vanishing and needs $X$ to be Calabi--Yau.   Moreover Zariski decomposition type arguments then imply that for some effective real divisor $\calE$ supported on  rigid non-movable surfaces, the real class $\Delta = D - 
\calE$ lies in the movable cone, and hence the Hessian is non-negative there.  This part of the argument may be extended to the general case of the \CY threefold $X$ having $b_2 (X) \ge 3$.

(3) In the case of $b_2 (X) =3$, our knowledge of the elliptic curve shows that $\Delta$ lies in a component of the positive index cone; here we use again that $E_i^3 \le 0$ for all 
the rigid non-movable surfaces (or more precisely that the $E_i$ lie in the  extended upper half-space).
 In particular the movable divisor $\Delta$ has $\Delta ^3 \ge 0$, and then a connectedness argument says that it lies in the \it same \rm component of the positive index cone as the hybrid one containing the \Kahler cone.
 We note however that the Hessian is negative at $D - \varepsilon \calE$ for $\varepsilon$ small and positive; this latter fact follows since $D$ is not visible from any of the $E_i$. 
 This gives a contradiction to the convexity of the component since $\varepsilon \Delta + (1-\varepsilon ) D$ is not in this component for 
 $0<\varepsilon  \ll 1$.

\smallskip
For $\rho =3$ the cubic defines a curve in the real projective plane, and our assumptions say that this is a real elliptic curve.  To study real 
elliptic curves, the Hesse normal form for the curve will be useful,  the theory of which may be found in 
Theorem 6.3 \cite{BM} or Section 3 of \cite{Dolg}.  Normally we might take 
real coordinates so that the real elliptic curve takes the classical
Hesse normal form  
$x^3 + y^3 + z^3 = 3kxyz$
with parameter $k$, but for our purposes it will be more convenient to make a change of coordinates, 
closely related to the canonical form described in Remark 6.9 of \cite{BM}, 
 so that the cubic takes the form 
$$ F(x,y,z) = -x^3 - y^3  - (z-x-y)^3 + 3kxy(z-x-y), \eqno{(1)}$$
and hence the `triangle of reference' of $F$ is now in the affine plane $z=1$ with vertices $ (0,0), (1,0)$ and $ (0,1)$.  We shall write $F_k$ if we wish to indicate the dependence on $k$.  
 An easy check verifies that the cubic $F$ is strictly negative at all points on the three affine lines $x=0$, $y=0$ and $x+y=1$. 
Recall that if $k>1$, then the real curve $F=0$ has two components, the bounded component (lying in the triangle of reference given by the above three real lines) and the unbounded component.  The cone 
in $\R^3$ corresponding to the bounded component has two connected components when one removes the origin, a positive part inside which $F>0$ and a negative part inside which $F<0$, whilst the cone on the unbounded component only has one 
connected component in $\R^3$, even after removing the origin.
 In the case $k>1$, the unbounded component 
 has three affine branches, one of which  
lies in the negative quadrant $x<0,\  y<0$, one in the sector $ y<0,\  x+y >1$ and the third in the sector $x<0,\   x+y >1$.  The 
(real) inflexion points of the cubic are at 
$B_1 = (0:1:0)$, $B_2 =(1:0:0)$ and $B_3 = (1:-1:0)$, i.e. the intersection of the line at infinity $z=0$ with the curve (a further reason why the chosen change of coordinates is helpful).  The asymptotes for the affine branches of the unbounded component may be found by calculating the tangents to the curve at the inflexion points, and are  
$$ x = - {1\over {k-1}}, \quad  y = - {1\over {k-1}} \quad \hbox{\rm and}\quad x+y = {k\over {k-1}}. \eqno{(2)}$$

Noting that $x^3 + y^3 + z^3 -3xyz = (x+y+z)(x+\omega y + \omega^2 z)(x+ \omega ^2 y + \omega z)$, where $\omega$ is a primitive cube root of unity, 
 when $k=1$ the cubic (1) splits into the real line $z=0$ and two complex lines 
(meeting at the centroid $({1\over 3}:{1\over 3}:1)$ of the triangle of reference).

When $k <1$, the cubic $F=0$ is smooth but with only one real component, with three affine branches, one in the region $x>0,\ y>0, \ x+y >1$, one in the region 
$x>0, \ y<0, \ x+y <1$ and one in the region $x<0, \ y>0, \ x+y < 1$.  The asymptotes are calculated as before and are given by the equations (2).  

 \begin{figure}
 \centering
     \includegraphics[width=10cm]{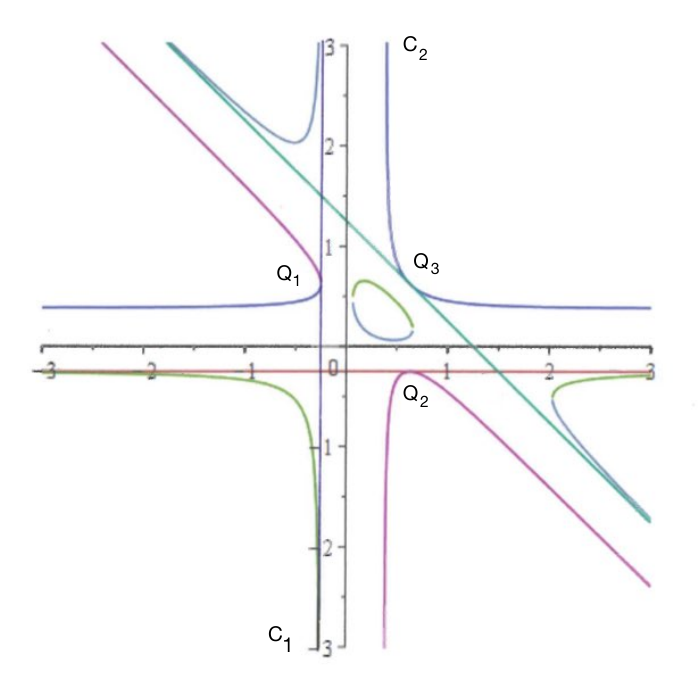}
\caption{Cubic $F_k =0$ with asymptotes and Hessian $H_k$ for $k=5$}
%\label{}
\end{figure}

By Remark 2.11 of \cite{BM},  the Hessian of the cubic $-x^3 -y^3 -z^3 + 3kxyz$ is given by  
$$27\bigl(2k^2(x^3 +y^3 +z^3 )- ( 8-2k^3) xyz\bigr).$$
Thus if $H_k$ denotes the Hessian of the cubic $F_k$, the fact that our change of coordinates was unimodular shows for $k\ne 0$ that $H_k = -54k^2F_{k'}$, with parameter $k' = {{4-k^3}\over {3k^2}}$.  In particular we see that if $k>1$, then $k' <1$, and so the Hessian curve of a real elliptic curve with two components has only one component.  For $k<1$, we have two notable values: $k=0$ for which the Hessian curve is the three real 
lines given by $-xy(z-x-y) =0$, and $k=-2$ for which the Hessian curve is the three lines 
(two of them complex) corresponding to $k' =1$ described before.  Apart from these two values, for any real elliptic curve with one component, the Hessian curve is a real
elliptic curve with two components, the bounded component lying in the triangle of reference.  
We shall assume for simplicity in this paper that the Hessian cubic is irreducible, namely that $k \ne -2$ or $0$; note that the above calculations then imply that $H$ is 
is strictly positive on the three affine lines $x=0$, $y=0$ and $x+y =1$.

The case of the curve having two components is illustrated in Figure 1 (which shows $F_5$, $H_5$ and the
three asymptotes for $F_5$).  Here the Hessian $H_5$ consists of the three branches (including $C_2$) containing respectively $Q_1, Q_2, Q_3$ and the real cubic $F_5 =0$ 
consists of the other four branches illustrated, namely the bounded component and three other unbounded components, including $C_1$. 
The picture for one component when $k \ne -2$ or $0$ is not dissimilar, where the roles of the cubic and its Hessian are interchanged --- see 
for instance Figure 3 below.  
    For $k \ne -2$, $0 $ or $1$, we note that the asymptotes to the cubic $F_k =0$ are tangent to the affine Hessian curve $H_k =0$; this is just a special case of a classical result that the 
double polar with respect to the cubic at a point on its Hessian is tangent to the Hessian at the image of the point under the Steinerian involution 
(see \cite{Dolg}, Section 3.2 and Exercise 3.8) --- the Steinerian involution on the Hessian will be explained at the start of Section 2.  When the point is 
an inflexion point of the cubic, the double polar is just the tangent line to the cubic, in our case the asymptote --- the corresponding affine points on the Hessian 
where these lines are tangent are labelled $Q_i$, $i=1,2,3$.  This gives more precise information about the 
affine regions where the Hessian curve can lie.

We can now identify the components of the positive index cone when $b_2 (X)=3$; this would be much more complicated for $b_2(X) >3$.  Taking $A= (a,b,1)$ in 
 the affine plane $z=1$, we wish to know the index of the quadratic form defined by the homogeneous quadratic  $G_A(D) = A\cdot D^2$; 
if $D = (x,y,z)$ in the above coordinates, this is 
explicitly given by 
 $$-ax^2 -b y^2 - (1-a-b)(z-x-y)^2 +kay(z-x-y) +kbx(z-x-y) + k(1-a-b)xy. $$

If $F_k =0$ has two real components, i.e. $k>1$, then this quadratic form at $A$ has index $(1,2)$ if
$H(A) >0$, index $(2,1)$ if $H(A) <0$ and index $(1,1)$ if $H(A)=0$.  
This may be easily verified by considering sample points such  as $A= (0,0,1)$, where the index is plainly $(1,2)$ 
and points $A = (t,t,1)$ for $t\gg0$ where the index is $(2,1)$.
One is therefore looking for the regions for which either $F>0$ and $H>0$, or $F<0$ and $H<0$, the latter being relevant since for $D$ in such a region, both $F$ and $H$ will be positive at $-D$.    

The  cubic curve then bounds precisely four (convex) regions of the (real) affine plane $z=1$ on which both $F>0$ and $H>0$ 
 including the bounded  component inside which $F>0$ (contained in the triangle of reference), where the index 
of the corresponding quadratic form is $(1,2)$.  Moreover the Hessian curve bounds precisely three (convex) regions of the affine plane on which both $F$ and $H$ are 
negative (where the index of the associated quadratic form is $(2,1)$).  Each unbounded affine  region on which 
 $F>0$ will together with the negative of the appropriate affine region 
on which  $H<0$ give rise to a connected component of the positive index cone in $\R^3$, part of whose boundary is contained in  $F=0$ and part of whose boundary is contained in  $H=0$, with the two parts meeting along rays corresponding to 
 two of the inflexion points of the curve $F=0$.  Moreover such a cone is also convex, since the tangent line to the 
 cubic at $B_i$ only meets the Hessian again at the point $Q_i$ (on the affine branch not containing $B_i$). 
 For each of these three resulting \it hybrid \rm  cones in $\R^3$, we have  $F>0$, $H>0$ and a continuity argument verifies that the index 
 is $(1,2)$ on each cone.  The other 
component of the positive index cone corresponds to the bounded component and whose boundary is 
contained in $F=0$.

In the case when $F=0$ only has one connected component, i.e. $k<1$, we have that the 
condition $F>0$ defines three (convex) regions of the affine plane $z=1$, 
and in these regions the Hessian is positive and 
the index is $(1,2)$.  Under the assumption that the Hessian is also irreducible,  
 there are four (convex) regions of the affine plane on which the 
Hessian is negative (where the cubic is also negative),
 three unbounded regions on which the index is $(2,1)$, and the region determined by the bounded component 
 of the Hessian, on which $F<0$ and $H<0$.  
  As before we obtain three \it hybrid \rm components of the positive index cone,
 obtained from unbounded affine regions on which $F>0,\ H>0$ together 
with the (negative of) unbounded affine regions on which $F<0,\ H<0$; a continuity argument again ensures that the index is $(1,2)$ on the hybrid components.

For  the Fermat cubic, i.e. $k=0$, we check easily that the index at any interior point of the triangle of reference is in fact $(0,3)$.  Thus by a continuity argument, for any $-2 < k < 1$, the index of the quadratic form for any $A$ inside the bounded component of the Hessian remains $(0,3)$, since as $k\to 0$ the bounded component of the Hessian of $F_k$ tends to the triangle of reference.  For points actually on this bounded component of 
the Hessian the index is $(0,2)$, apart from the Fermat cubic and the vertices of the triangle of reference, where 
the index is $(0,1)$.  
  For $k<-2$, we check that for $A$ inside the bounded component, the index 
is $(2, 1)$ (check it at a suitable such point for $k \ll -2$ and use continuity again), and for $A$ 
actually on the bounded component the index is $(1,1)$.  Thus in this case, 
 points in the negative of the corresponding open cone have index $(1,2)$, and we obtain  a further connected component of the positive index cone, whose boundary is contained in $H=0$.

In both cases under consideration, namely the real elliptic curve has two components or one component with irreducible Hessian, the components of the positive index cone are all convex, and their closures are strictly convex cones.

Given a real smooth ternary cubic, 
we suppose that we have fixed the real coordinates chosen in equation (1) to transform the cubic into the given Hesse form.  This then determines an open upper half-space 
given in coordinates by real classes $E = (a,b,c)$ with $c>0$.  This half-space may be defined in a coordinate-free way, since the three real inflexion points of the elliptic curve determine a real projective line,  and the corresponding plane divides $\R^3$ into two open half-spaces,
the upper half-space identified as where the cubic being positive defines  convex components.

\begin{defn}  Given a smooth ternary real cubic, the \it extended upper half-space \rm with respect to the cubic is the open upper half-space just identified together with those boundary points representing an inflexion point of the elliptic curve

\end{defn}

\begin{lem} Given a smooth ternary real cubic with irreducible Hessian, 
and a real class $E$ at which the index is $(1,q)$ for $q\le 2$  (for instance in the case which interests us,  the class 
of a rigid non-movable surface) with $E^3 \le 0$, then $E$ lies in the extended upper half-space.

\begin{proof} We suppose that coordinates have been chosen as in equation (1), and $E= (a,b,c)$ with respect to these coordinates.
If $c\le 0$ and $E$ does not represent an inflexion point of the elliptic curve, then the fact that $H(E)\ge 0$ implies that $F(E) >0$.
\end{proof} \end{lem}

We shall prove below all but the last sentence of Theorem 0.2 under the weaker assumption that the rigid non-movable surface classes 
lie in the  extended upper half-space  --- the only place that we shall need the stronger property that $E^3 \le 0$ for all rigid non-movable surface classes is in proving Corollary 1.7, which in turn implies the last sentence of Theorem 0.2.    
Note that because we have only chosen a \it real \rm coordinate system, it is unclear when a point $(a,b,c)$ represents a class in $H^2 (X, \Z)$; this however will not be a problem in our arguments.

When $-2 < k < 0$ or $0<k<1$, we know from Hodge index considerations and the above calculations that the negative  closed convex 
cone on the bounded component of the Hessian $H$ cannot contain the class of any surface, and in particular 
the cone cannot contain the \Kahler cone of $X$.  Under our assumptions, 
we show below that this latter fact continues to be true when $k<-2$.  

Recall that for a convex body $V$ in $\R^n$ with non-empty interior, a point $D$ on the boundary $\partial V$ is said to be \it visible \rm (with respect to $V$) from a point $A$ if the line segment $AD$ does not meet the interior of $V$.
  If $D$ is a smooth point of $\partial V$, the tangent hyperplane 
through $D$ determines a closed half-space whose intersection with the interior of $V$ is empty;
 then $D$ is visible from $A$ if and only if $A$ is in this half-space.   If $W$ is a convex cone with vertex at the origin and 
 non-empty interior, then for any $A\ne 0$, the set of points 
in $\partial W$ which are visible from $A$ and $-A$ will be called 
the \it visible extremity \rm of $\partial W$ from $A$; a non-zero smooth point $D \in \partial W$ is in the visible 
extremity from $A\ne 0$ if and only if $A$ is in the tangent hyperplane.
 A non-zero $D \in  \partial W$ which is not in the visible extremity from  $A\ne 0$ will be visible from $-A$ if and only if it is not visible from $A$, and also if and only if  $-D \in \partial (-W)$ is visible from $A$ with respect to $-W$.

\begin{prop}   Assuming that the rigid non-movable surface classes lie in the  extended upper half-space, if the cubic form corresponds to the case $k<-2$, and so 
 there is a component $P^\circ$ of the positive index cone corresponding to the bounded component of the Hessian, then the \Kahler cone of $X$ is not contained in $P^\circ$.
 
 \begin{proof} We suppose that the \Kahler cone is contained in $P^\circ$.  
 Index considerations imply that the interior of the movable cone is then contained in $P^\circ$, and without loss of generality we may assume that $X$ is general in moduli.  If there are no rigid non-movable surfaces on $X$, then any effective divisor must be movable, as each component then corresponds to a nef (and hence mobile) divisor on some minimal model --- see the first paragraph of Section 2 from \cite{WilBd}. 
 Choosing any $D$ on the boundary of the movable cone,   by the second proof of Theorem 0.1 in Section 4 of \cite{WilBd},  we deduce $\vol (D) \ge D^3 >0$, and hence $D$ is big;
 this remains true for all nearby rational points $D'$, 
 including those not in the movable cone, therefore yielding a contradiction.  
 
 If there is precisely one rigid non-movable surface $E$ on $X$,  we choose any $D$ on the boundary of the movable cone which is not visible from $E$;  then the above argument shows that  
 $\vol (D')>0$ for all nearby rational points, including $D' = D - \epsilon E$ with $0<\epsilon \ll1$ 
 not in the movable cone.  Since $D'$ is big, we have that $D' - \lambda E$ is movable for some $\lambda > 0$, 
 and this is a contradiction, since $D'$ is a convex combination of $D' -\lambda E$ and $D$.  We comment in passing that this argument does not use the assumption 
 that $E$ lies in the upper half-space.
 
 Next we suppose that 
 there are precisely two rigid non-movable surfaces $E_1, E_2$ on $X$, both in the  extended upper half-space, 
 where $E_2$ is not allowed to be a negative multiple of $E_1$ in the case of points with $z=0$ (i.e. representing inflexion points).
 We observe that the closed  line segment joining $-E_1$ and $-E_2$ is disjoint from the movable cone, since 
 a convex combination of $E_1$ and $E_2$ is effective and so its negative cannot also  be pseudo-effective, and hence cannot lie in the movable cone.  We can then find a point $D$ on the boundary of the movable cone 
 which is visible from both  $-E_1$ and $-E_2$, in particular with the interior of the cone generated by $D$, $-E_1$ and $-E_2$  disjoint from the movable cone.  The previous argument implies that 
 $\vol (D') >0$ for all nearby rational points,  and in particular we may take a big such $D'$ in the interior  of the cone generated by $D$, $-E_1$ and $-E_2$.
 Writing such a $D'$ in terms of its movable part $\Delta$ and non-negative rational multiples of the $E_i$ yields a contradiction, since then 
  $\Delta$ could not lie in the movable cone.

 Finally we may suppose therefore that  
 there are at least three rigid non-movable surfaces on $X$, say $E_i$ for $i = 1,2, 3$,  and let $L\in P^\circ$ denote a very 
 ample class.  These four (effective) 
 integral classes are 
 linearly dependent in $H^2 (X, \Z)$, and we cannot have any $E_i$ being a rational convex combination of the other three classes since it is non-movable.  Also $L$ is not 
  a rational convex combination of the $E_i$ since they lie in different half-spaces.  We deduce without loss of generality that some integral convex combination of say $E_1$ and $E_2$ is an integral convex combination of $E_3$ and $L$, and hence is mobile.  As however the cone generated by mobile classes lies inside the closed cone $P$ by index considerations, and 
 hence lies in the open lower half-space $z<0$, this is a contradiction.  
   \end{proof}
 \end{prop}
We now consider the case where the real elliptic curve has two components, and the \Kahler cone is contained in the positive cone on the bounded component.  We shall need an elementary lemma in convexity 
theory, the idea of the proof given being suggested to the author by a colleague Imre Leader.

\begin{lem} Let $V \subset \R^2$ be a open bounded convex body with a smooth boundary curve, and $x_1 , x_2 ,  \ldots $ be a 
(perhaps infinite) collection of points in $\R^2$ not in $V$.  We let $W$ denote the points of
the boundary $\partial V$ which cannot be seen  from any of the $x_i$, i.e. $W$ consists of points $x \in \partial V$ such that the line segments $xx_i$ all meet $V$.
Then $V$ is contained in the closure $Z$ of the convex hull of $W$ and $x_1, x_2 , \ldots $.

\begin{proof} We suppose that the result is not true, and so in particular by the convexity of $V$ there is a point $x$ in the boundary of $V$ which is not in the closure $Z$ of the convex hull of $W$ and $x_1, x_2 , \ldots$.  Standard convexity results imply the existence 
of a  line $l$  though $x$ such that $Z$ is strictly on one side of $l$, and in particular all points of $W$ and all the $x_i$ are in the corresponding open half-plane.  

If $l$ is tangent to the boundary $\partial V$ at $x$, then $V$ itself must be on one side of $l$.  If $W$ and 
the $x_i$  are contained in the open half-plane disjoint from $V$, we get an immediate contradiction by 
considering $y$ the other point on $\partial V$ whose tangent is parallel to $l$, which 
therefore cannot be seen from any of the $x_i$ and hence by definition lies in $W$, contrary to assumption.  If however 
$W$ and the $x_i$  are contained in the same open half-plane that contains  $V$, then  $x$ cannot be seen from any of the $x_i$ and so $x\in W$, contrary to assumption.

If however $l$ is not tangent to the boundary $\partial V$ at $x$, we consider the point $y\in \partial V$ with tangent parallel to $l$ such that $y$ is on the other side of $l$ to $W$ and the $x_i$.  Then 
$y$ cannot be seen from any of the $x_i$ 
and hence $y\in W$, a contradiction.
\end{proof}
\end{lem}

\begin{prop}  Suppose that the cup-product cubic form is smooth with irreducible Hessian, and that all rigid non-movable surface classes lie in the  extended upper half-space.
Suppose also that $X$ is general in moduli and that the real elliptic curve has two components, with the \Kahler cone of $X$  contained in the open positive cone $P^\circ$ on the bounded component.  Then $P^\circ$ is contained in the 
interior of the effective cone of $X$.

\begin{proof}
Consider now the affine plane $z=1$ and let $V$ be the open convex body 
in the affine plane given by the bounded component of the above real elliptic curve, so that $P^\circ$ is the cone on $V$.  Suppose that $E = (a,b,c)$ represents the class of a rigid non-movable surface; by assumption, either $E$ represents an inflexion point of the elliptic curve  
or $c>0$ and there is a unique point $A$ of the affine plane which is a positive multiple of $E$;  
moreover $A\not\in V$  from Proposition 4.4 of \cite{WilBd}.  
The points of the boundary of $V$ which cannot be seen from $E$  with respect to $P^\circ$ are precisely those points $D$ with 
$E\cdot D ^2 >0$; the points at which $E\cdot D^2 =0$ corresponding to the tangent lines to the boundary which pass through $A$, with an appropriate interpretation for 
the limit case of a 
 point $E$ at infinity on the plane $z=0$ in $\R^3$ (therefore representing an inflexion point).  
We let $\tilde Q$ denote the (convex) subset of $V$ defined by the inequalities $E_i \cdot D^2 >0$ for all $i$, 
which contains 
in its closure the set of all $D_0 \in \partial V$ which cannot be seen 
from any $E_i$ 
(note that the cone $Q$ on $\tilde Q$ contains the \Kahler cone).  Given such a $D_0$, for any ample class $L$, note that any strictly convex combination of 
$D_0$ and $ L$ lies in $\tilde Q$.  This enables us to find rational points $D_i \in \tilde Q$ for $ i>0$ with $D_i \to D_0$.  By the argument in Proposition 4.1 and Lemma 4.3 from \cite{WilBd}, each $D_i$ is in the effective cone (the quoted results here  use the assumption that $X$ is general in moduli but do not need there to be only finitely many rigid non-movable surfaces), and so the point $D_0$ is in the pseudo-effective cone.  Thus, from Lemma 1.5, it follows by a limiting argument that the closure $P$ of $P^\circ$ is in the 
pseudo-effective cone, and hence the open cone $P^\circ$ is contained in the interior of the effective cone (i.e. the \it big \rm cone).
\end{proof}
\end{prop}

\begin{cor}  Under the hypotheses of Theorem 0.2, the linear form $c_2$ is non-negative on the positive cone $P$ on the bounded component of the elliptic curve.

\begin{proof} Without loss of generality, we may assume that $X$ is general in moduli.
We noted in Section 2 of \cite{WilBd} that if $c_2 \cdot E <0$ for some rigid non-movable surface, 
then $E^3 >0$.  Thus the hypotheses imply that $c_2 \cdot E \ge 0$ for all rigid non-movable surfaces $E$ on $X$.  Moreover we also noted 
in Remark 1.2 of \cite{WilBd} that $c_2 \cdot L \ge 0$ for any movable class $L$.  Since by the previous result any element of $D \in P^\circ$ is effective, and $X$ is general in moduli, 
it follows that $c_2 \cdot D \ge 0$ for any $D\in P$.
\end{proof}
\end{cor}
In the light of these results, for Theorem 0.2 it is sufficient to prove:

\begin{thm} Assuming that the cup-product cubic form is smooth with irreducible Hessian, and that all rigid non-movable surface classes lie in the  extended upper half-space, 
 the \Kahler cone is not contained in a hybrid component of the positive index cone.
\end{thm}

We note that in the explicit example given in Remark 0.1, it is straightforward to check that the Hessian matrix is negative definite at say $(-4,2,1)$, from which we 
may deduce that the invariant $-2 < k <1$ and the \Kahler cone must be contained in a hybrid component of the positive index cone.

For the examples of 6-manifolds not supporting any \CY structure, we use Wall's Theorem; 
we may take the cubic form on $H^2 (M, \Z ) = \Z^3$ to be an appropriate 
integral multiple of $F$ as in equation (1), with $k<1$ an integer other than $-2$ or $0$, and any suitable integral linear form satisfying the congruence conditions --- we note that the trilinear form corresponding to $F$ does satisfy the first of the congruence conditions.  We may also 
take the cubic form on $H^2 (M, \Z ) = \Z^3$ to be an appropriate integral  multiple of $F$ as in equation (1) with 
any integer $k>1$ and any suitable integral linear form satisfying the congruence conditions but which is negative somewhere 
on the open positive cone on the bounded component of the real elliptic curve.  In these cases, if the 
integral multiple of $F$ 
has been chosen so that there are no integral classes $E$ with $1\le E^3 \le 8$, or indeed as we saw in the Introduction if we scale both $F$ and a linear form satisfying the required congruence conditions by say a factor of 4, 
our Theorem 0.2 rules out 
the possibility of any \CY structures, independent of any choice for the 
even rank free abelian group $H^3 (M, \Z )$.   Being even more concrete, the integral cubic  form $nF_k$ for any integer $k<1$ ($k\ne -2$ or 0) and any integer $n>8$ 
can for instance never be the cubic form arising from a  \CY threefold.

\begin{notation}  We shall now fix on the notation that will be used in the rest of this paper to describe a hybrid component $P^\circ$ of the positive index cone, with $P$ denoting its closure.  
Taking the affine slice $z=1$, where the cubic is of the form described above, with slight abuse of notation concerning the points at infinity, 
we denote the  projectivised boundary of $P$ by $C = C_1 \cup C_2$, where $C_1$ is an 
affine  branch of $F=0$ and $C_2$ is an associated affine branch of $H=0$, 
with $C_1$ and $C_2$ meeting at two inflexion points (at infinity).  When $k>1$, without loss of generality we may by symmetry take $C_1$ in the negative quadrant and then $C_2$ in the region $x>0, \ y>0, \ x+y >1$. 
The branches $C_1$ and $C_2$ meet (at infinity) at the 
inflexion points $B_1 = (0:1:0)$ and $B_2 = (1:0:0)$.  More specifically, the cone $P$ then has boundary the positive cone on $C_1$ together with the negative cone on $C_2$, the two parts meeting in  two  rays corresponding to positive multiples of $(0, -1, 0)$ and $(-1,0,0$).  When $k<1$ and $k \ne 0, -2$, we take $C_1$ in the region $x>0,\ y>0,\ x+y >1$ and 
$C_2$ in the negative quadrant. The branches $C_1$ and $C_2$ again meet (at infinity) at the 
inflexion points $B_1$ and $B_2 $.  The cone $P$ then has boundary the positive cone on $C_1$ together with the negative cone on $C_2$, the two parts meeting in the two  rays corresponding to positive multiples of $(0, 1, 0)$ and $(1,0,0$).  In both cases, we shall denote by $V_1$ the open convex subset of the affine plane bounded by $C_1$ and $V_2$ the open convex subset bounded by $C_2$. 
\end{notation}

\begin{prop}  Assuming that the cup-product cubic form is smooth with irreducible Hessian, and that all rigid non-movable surface classes 
lie in the  extended upper half-space, suppose that $P^\circ$ is a hybrid component of the positive index cone as described in Definition 1.1, $E_1 , E_2 , \ldots $ are the (perhaps infinitely many)  rigid non-movable classes in $H^2 (X, \Z )$ and $Q$ 
a connected component of the subcone of $P^\circ$ defined by the inequalities ${E_i \cdot D^2 >0}$ for all $i$.  
If $Q$ contains the \Kahler cone
then there cannot exist a non-trivial open arc of points $-D \in C_2$  which are visible with respect to the cone  $-Q$ from every $E_i$, but with each $D$ representing a point of the boundary $\partial Q$.
\begin{proof}  Again, we may without loss of generality assume that $X$ is general in moduli.  
We have that $Q = \bigcap _{i\ge 1} Q(i)$, where $Q(i)$ is the component of the subcone of $P^\circ$ defined by $E_i \cdot D^2 >0$ which contains the \Kahler cone.  
Were such an arc of points in $C_2$ to exist, we choose a point $-D_0$ in this arc; note that no point 
$E_i$ lies on the tangent plane to the cone on $C_2$ along the ray $\R _+ (- D_0)$, since otherwise some points of the arc would not be visible from $E_i$.  
 Thus $D_0$ is not visible from any of the $E_i$ with respect to $Q$.
For any real ample divisor $L$, we note that any strictly convex combination of $D_0$ and $L$ lies in each $Q(i)$, and hence  in $Q$.   In this way we show, using 
the argument in Proposition 4.1 and Lemma 4.3 from \cite{WilBd} 
as in the proof of Proposition 1.6, where the quoted results use that $X$ is general in moduli, 
that $D_0$ is a limit of effective rational divisors $D_j \in Q$; therefore $D_0$ is pseudo-effective. 

For any prime divisor $\Gamma$, we have a function $\sigma_\Gamma$ on the pseudo-effective cone as 
in Definition 1.6 of Chapter III from \cite{Nak} --- all the following references to \cite{Nak} will be from Chapter III.
Moreover by Proposition 1.10 and Corollary 1.11 of the given Chapter, there are at most three $\Gamma$ 
with $\sigma _\Gamma (D_0)>0$, whose classes are moreover linearly independent 
in $H^2 (X, \R )$.  Under our assumption that $X$ is general in moduli these will be rigid non-movable surfaces $E_i$, 
without loss of generality 
 $E_1 , \ldots ,E_r$ with $r\le 3$.  Thus in the $\sigma$-{\it decomposition }\rm (see Definition 1.12 from the cited Chapter)
$D_0 = P + N$ of $D_0$, we have $N = \Sigma _{i=1}^r 
\sigma _{E_i} (D_0) E_i$.  In particular, by Lemma 1.8 of the Chapter, we have $\sigma _{E_i} (P) = 0$ for 
$i=1,\ldots ,r$, and indeed also that $\sigma _{\Gamma} (P) = 0$ for any other surface  $\Gamma$.
If we knew that $D_0$ is big, then by Lemma 1.4 (4) of the cited Chapter,
applied $r$ times, we have that $D_0 = \Delta +  \Sigma _{i=1}^r 
\sigma _{E_i} (D_0) E_i$, with the real divisor $\Delta$ also big and in particular pseudo-effective, and 
$\sigma _\Gamma (\Delta) = 0$ for all surfaces $\Gamma$.  Thus by Lemma 1.14 (1) of the Chapter, we 
would deduce that $\Delta$ is movable.

Suppose first that $r=0$; then by Lemma 1.14 (1) of \cite{Nak} Chapter III, we deduce that $D_0 $ is movable.  We saw
then  in the second proof of Theorem 0.1 in  Section 4 of \cite{WilBd} 
that $\vol (D_0) \ge D_0^3 >0$ and hence $D_0$ is big.  If there are no rigid non-movable surfaces $E$ on $X$, we noted in \cite{WilBd} that this gives an immediate contradiction, since then there would exist (rational) points $D$ near $D_0$ at which the Hessian is negative but which are big and hence (recalling that $X$ is assumed general in moduli) in this case also movable;  
this is a contradiction.  However, 
 if there exists a rigid non-movable surface $E$, we deduce that $D_0  -\epsilon E$ is also big for $0<\epsilon \ll1$, and hence (using a previous argument) of the form $\Delta + \calE'$, with $\Delta$ big and movable and $\calE'$ supported on at most three of the $E_i$. Therefore 
$D_0  = \Delta + \calE' + \epsilon E$.

We suppose now that $r>0$ and  
verify that $D_0$ is big. 
We know that there is an $E= E_i$ with $\sigma _E (D_0)>0$; choose $\delta = \sigma _E (D_0)/2$.  
Thus, given a real ample divisor $L$, for all $0<\epsilon \ll 1$ the big divisor $D_0 + \epsilon L$ has $\sigma _E (D_0 + \epsilon L)> \delta$, straight from the definition of $\sigma _E$ 
as a limit in Definition 1.6 of Chapter III 
from \cite{Nak}.  Then by Lemma 1.4 (4) of the Chapter, 
 we deduce that $D_0 + \epsilon L -\delta E $ is big for all $0 <\epsilon \ll1$.  Thus $D_0 -\delta E$ is pseudo-effective.  From our choice of $D_0$, we know that $D_0$ is not visible from $E$ with respect to $Q$ and hence 
$D_0 +tE \in Q$ for $0<t \ll 1$, and hence effective (as we saw above).  Since \it all \rm points $-D$ of the original arc in $C_2$ have $D$ pseudo-effective, we must have that $D_0$ lies in the interior of the pseudo-effective cone, and hence is big as claimed.  We saw above that then
 $D_0 = \Delta + \calE$, for some big class $\Delta$ in 
the movable cone and some real non-zero 
 convex combination $\calE$ of at most three of the $E_i$.

Summing up therefore, we can in all cases  
write $D_0 = \Delta + \calE$, for some big class $\Delta$ in 
the movable cone and some real non-zero 
 convex combination $\calE$ of finitely many of the $E_i$, which by assumption lies in the upper half-space. 
 Thus $\Delta = D_0 - \calE$ lies in the open  lower  half-space; since the Hessian at $\Delta$ is non-negative, it follows that 
 $\Delta ^3 >0$ (cf. Lemma 1.3).  If $L$ denotes an ample divisor in $P^\circ$, then $\Delta +tL$ is movable for all $t\ge 0$, and 
 $(\Delta +tL)^3 > 0$ for all $t\ge 0$, from which it follows by connectedness that $\Delta = D_0 - \calE \in P$. 
 Recall that $C_2$ has been assumed smooth; since $-D_0 \in C_2$ is visible from every $E_i$ with 
respect to $V_2$, 
it follows that the Hessian at $-D_0 + \varepsilon \calE$ is strictly positive for $0< \varepsilon \ll 1$, and hence 
the Hessian at $D_0 - \varepsilon \calE$ is strictly negative; this contradicts the convexity of  $P$. 
\end{proof}
\end{prop}

Our technique for ruling out such a hybrid component $P^\circ$ containing the \Kahler cone is as follows; 
in the next three sections of the paper, we prove the following result in the various cases for the elliptic curve.
In the light of this combined with Proposition 1.9, the only conclusion then is that if the rigid non-movable classes all lie in the 
 extended upper half-space, then
  $P^\circ$ cannot contain the \Kahler cone, and so Theorem 1.8 will be proven.  The proof of Proposition 1.10 is slightly 
technical when written out in full, but the basic mathematical  input consists just of 
classical facts about the Steinerian involution on the Hessian of a real elliptic curve; the basic  properties of the $E_i$ we use are that the index at $E_i$ is $(1,q_i)$ with $q_i \le 2$, and that either $E_i$ represents an inflexion point of the cubic or it lies in the open upper half-space.

 \begin{prop} Assuming that the cup-product cubic form is smooth with irreducible Hessian, and that every rigid non-movable surface class lies in the  extended upper half-space,
 let $P^\circ$ be a hybrid component of the positive index cone.  Suppose  $E_1, E_2, \ldots $ denote 
 the (perhaps infinitely many) rigid non-movable surface classes in $H^2(X, \R)$  and let $Q$ 
 be  a connected component of the subcone of $P^\circ$ defined by the inequalities $E_i\cdot D^2 >0$ for all $i$,
 where $Q$ is assumed to have  non-empty interior.  Then there exists a non-trivial open arc in $C_2$ of points $-D$ which are visible from all the $E_i$
 with respect to the cone  $-Q$, 
 with each $D$ representing a point of the boundary $\partial Q$.
\end{prop}

%\begin{rem} In the special case where there is at most one rigid non-movable surface $E$ on $X$, using the methods of this paper we can show that the previous result remains true irrespective of the value of $E^3$.  Thus Theorem 1.8 and 
%(in view of a comment in the proof of Proposition 1.6) all but possibly the last sentence of  Theorem 0.2 hold true in this case also.
% \end{rem}

  Theorem 1.8 states that if the \Kahler cone is contained in a hybrid component, then there must exist rigid non-movable surfaces $E$ with $E^3 >0$ lying in the closed lower half-space.  There will however exist at most two of these.
  
  \begin{prop}  Suppose $X$ is a \CY threefold with $b_2 (X)=3$ and 
  whose corresponding cubic is smooth with irreducible Hessian, with the \Kahler cone  contained in a hybrid component of the positive index cone.  
  Then there are either one or two rigid non-movable surfaces $E$ with $E^3>0$ lying in the closed lower half-space.

  \begin{proof}
 Assuming that the K\"ahler cone is contained in the hybrid component $P^\circ$ 
(with the usual convention that the boundary is determined by the two affine curves $C_1$ and $C_2$), we know from Proposition 4.4 of \cite{WilBd} that for  any such class $E$, we have  $E\not \in P$.  Hence 
if $k>1$  any such class lies in the open half-space $x + y > \frac{k'}{k'-1}z$ (where as usual $k' = \frac {4-k^3}{3k^2}$) corresponding to the relevant asymptote to the Hessian.
If $k<1$ we have a tangent line $x+y =e_2$ to the bounded component of the Hessian with the bounded component on the opposite side to $C_2$ (see Section 3) and 
then any such class lies in the half-space $x+y \le e_2  z$, in the case of $k< -2$ possibly lying in the negative of the closed cone on the bounded component of the Hessian.  Suppose now that there were three classes $E_1, E_2, E_3$ as above; choose a very ample class $L \in P^\circ$;
since there is an integral dependence relation between these four (effective) classes, we deduce that some class $\calE$ which is a convex combination of two or three of the $E_i$ 
(and hence in the appropriate half-space  $x+y > \frac{k'}{k'-1} z$, respectively $x+y \le e_2  z$, for respectively $k>1$ and $k<1$) is 
in fact mobile (cf. the last part of the proof of Proposition 1.4), and hence has non-negative value for the Hessian.  We deduce that $-\calE$ is also  in the closed upper half-space and has non-positive Hessian.  If $k<-2$, we note that $\calE$ is not in the negative cone on the bounded component of the Hessian, since then by index considerations this cone would also 
contain the mobile cone. 
 Thus $-\calE$ lies in the closed cone on one of the two affine regions bounded by the unbounded branches of the Hessian other than $C_2$  (cf. Sections 3 and 4) and so 
 $\calE$ lies in a component of the positive index cone distinct from $P$.  
 This however is impossible: since $\calE^3 > 0$, any convex combination of $\calE$ and $L$ will be both mobile and have strictly 
 positive cube,   which by connectedness of $P$ would imply that $\calE \in P$ and hence a contradiction.  \end{proof}
\end{prop}

\begin {rem} This last result may be compared with the more general fact that if the cubic is smooth with irreducible Hessian, then $X$ contains at most six rigid non-movable surfaces $E$ with $E^3>0$.  To see this, observe that there are at most four components of the positive index cone and that the component containing the \Kahler cone cannot contain any of these surfaces by Proposition 4.4 of \cite{WilBd}.  Moreover, each $E$ is contained in the positive index cone and any of the 
other components of the positive index cone will contain at most two of these surfaces: if a 
component (disjoint from the \Kahler cone) contained three, say $E_1, E_2, E_3$, then by considering the linear dependence relation between these three classes and a very ample class $L$, we deduce that an integral combination $\calE$ of two of the $E_i$ is a \it mobile \rm class in the given component, and in particular has $\calE^3 >0$.  The line segment $\calE+tL$ then consists of 
mobile divisors with strictly positive cube, contradicting the fact that $\calE$ and $L$ lie in different components of the positive index cone.
\end{rem}

\section{Hybrid components when elliptic curve has two real components}

In this section, we study the case when $k>1$, i.e. the real elliptic curve $F=0$ has two components, and so in particular the Hessian is smooth.  For simplicity, from now on we often use $H$ to denote both the Hessian and the Hessian curve.

With notation as in the previous section, we let  $P^\circ$ denote a hybrid component of the positive index cone, with closure $P$.  Without loss of generality, we may by symmetry adopt the notation explained in the previous section, 
where the component is determined by the branches $C_1$ and $C_2$ in the affine plane given by $z=1$.  Recall also that we denote by $V_i$ (for $i=1,2$) the 
open convex affine region 
bounded by $C_i$.

Let $E$ denote a class of a rigid non-movable surface on $X$; thus the index at $E$ is $(1,q)$ with $q\le 2$ and by assumption  $E=(a,b,c)$ with $c\ge 0$, 
and  $c=0$ only if $E$ represents one of the inflexion points of the cubic.  We now specify where the homogeneous quadratic given by  
$G_E (D) = E\cdot D^2$  vanishes on 
 the boundary of  $P$. If $E \in \R^3$ represents one of the inflexion points, we see below that the function $G_E$ is  easy to understand explicitly.
  The main case to consider  is when $c>0$ and so some positive real multiple $A$ of $E$ lies in the affine plane $z=1$.  We  
 shall therefore need  to study real classes $A$ in the affine plane which have 
 index $(1,q)$ for $q\le 2$.   
 
 For $A$ in the affine plane, there is  a simple answer to the question 
of where on  $C_1$ we have vanishing of $G_A$, namely points on  $C_1$ for which the tangent passes through $A$ (including maybe inflexion points at infinity), noting that if $A\in V_1$, 
then $G_A$ is strictly positive on $C_1$.
There will be two such points if $A\not \in V_1$ is in the quadrant 
$x\le -{1\over {k-1}},\  y\le -{1\over {k-1}}$ (including the possibilities of the inflexion points $B_1$ and $B_2$, or a point on $C_1$ taken twice), 
no such points if $A$ is in the quadrant 
$x> -{1\over {k-1}}, \ y> -{1\over {k-1}}$, and one point otherwise.

We now recall some very classical theory:  when the Hessian curve $H$ 
of an elliptic curve $F=0$ is smooth, 
 there is a 
well-defined base-point free involution $\alpha$ on $H$, known as the  \it Steinerian map \rm or \it Steinerian involution \rm  (\cite{Dolg}, Section 3.2, noting a misprint in Corollary 3.2.5), where the polar conic of $F$ with respect to a point $U$ on $H$ is a line pair with singularity at
$U' =\alpha (U)$.  We note that this says that $U\cdot U' \equiv 0$ and $\alpha$ induces an involution on the real points of $H$.  The Steinerian involution  has the property 
that for any point $U' \in H$, the second polar of $F$  with respect to 
the point $U'$ is the tangent to the Hessian $H$ at $U= \alpha (U')$ (\cite{Dolg}, Exercise 3.8, again noting a misprint).
So for any 
  point $A$ on the tangent line to the Hessian  at $U$,  the conic $G_A =0$ contains the point $U'$; moreover if $A\ne U$, 
  the conic is non-singular at $U'$ with tangent line $L$ at $U'$ independent of $A$; explicitly $L$ is defined by the linear form $W\cdot U'$, where $W$ is any point  ($\ne U$) on the tangent line through $U$.  Moreover since such a linear form $W\cdot U'$ vanishes at both $U$ and $U'$, the common tangent line $L$ is just the line joining $U$ and $U'$.  
  Moreover $U'$ is clearly the unique point of intersection of the conic with $L$.

In the case currently under consideration, where the Hessian has only one real component, 
 the branch $C_2$ of the Hessian $H=0$  is in the region  $x>0,\  y>0, \  x+y >1$, passing through $B_1, B_2$ and $R = Q_3$, the affine point 
with coordinates $( {k\over {2(k-1)}}, {k\over {2(k-1)}})$.  
As $B_3 = (-1:1:0)$ is the third inflexion point, it is a consequence of the properties of the Steinerian map that 
the tangent (namely the double polar)  to $F=0$ at $B_3$ is the tangent to $H=0$ at $R = \alpha (B_3)$.  If therefore we take $B_3$ to be the zero of the group law, then $R$ is 
the unique real 2-torsion point of the Hessian and  
 $\alpha$ is given by translation in the group law by this point.    Let $Q_1$ be the point on the branch of $H$ in the region $x<0, \ y>0, \ x+y <1$, which is the 2-torsion point when we take $B_1$ as the zero in the group law, and  $Q_2$ the point on the branch of $H$ in the region $x>0, \ y<0, \ x+y <1$ 
corresponding to taking $B_2$ to be the zero is the group law.  We note that $\alpha (B_1) = Q_1$, $\alpha (B_2) =Q_2$ and $\alpha (R) = B_3$.  Thus the second polar of $F$ with 
respect to the inflexion point $B_1$ is the tangent to the Hessian at $Q_1$, namely given affinely by $x= -{1\over {k-1}}$, the 
second polar of $F$ with 
respect to the inflexion point $B_2$ is the tangent to the Hessian at $Q_2$, namely given affinely by $y= - {1\over {k-1}}$,  and the second polar of $F$ with 
respect to the point $R$ is the tangent to the Hessian at $B_3$, namely the asymptote $x+y = {k'\over {(k'-1)}}$, where as before $k' = {{4-k^3}\over {3k^2}}$.
Under the Steinerian involution, the arc $Q_1 B_3$ of the Hessian corresponds to the arc $B_1 R$ of $C_2$, whilst the arc $B_3 Q_2$  corresponds to the arc $RB_2$ of $C_2$.

Setting $A =(a,b,1)$, it is easily checked from this that $G_A (B_1) >0 $ if and only if $a< -{1\over {k-1}}$, that $G_A (B_2) >0$ if and only if  $b < -{1\over {k-1}}$ and 
that $G_A  (R)>0$ if and only if $a+b < {{k'}\over {(k'-1)}}$.  Moreover if   
$A=Q_i$ (for $i=1,2$)  then $G_A (B_i) =0$, and if $A= B_3$ then $G_A(R)=0$. We are interested in the cases of $B_1$ and $B_2$ 
since we want to know the sign of $G_A$ on points of $C_2$ where either $y \gg 0$ or $x\gg 0$.  
For future use, we  introduce the notation that  $\tilde B_1 = (0,1,0)$, $\tilde B_2 =(1,0,0)$ and 
$\tilde B_3 =(-1,1,0)$, points in $\R^3$ representing the three inflexion points on the real projective cubic.  

When $A$ is in the quadrant $x\le -{1\over {k-1}},\  y\le -{1\over {k-1}}$, then no tangent line at a point on the open arcs $Q_1B_3$ or $B_3Q_2$ contains $A$, and so 
$G_A$ is non-vanishing 
(and indeed positive)  on the affine branch $C_2$ --- in the case of equality perhaps vanishing at $B_1$ or $B_2$.  In this case we note that all of $C_2$ is visible from $A$ with 
respect to $V_2$.  
Moreover, when $E$ is a positive multiple of $-\tilde B_1$ or $-\tilde B_2$, then $G_E$ is positive on 
both $C_2$ and $P^\circ$, and all  of $C_2$ remains visible from $E$ with respect to $-P^\circ$.
When however $A$ is 
in the open region given by  $a> -{1\over {k-1}}$, $b > -{1\over {k-1}}$ 
and $a+b > {k'\over {k'-1}}$, then $G_A$ is non-vanishing 
on all the projectivised boundary of $P$, and hence $G_A$ is strictly negative on all of $P$ by continuity, since it plainly is for $A\in V_2$ by Lemma 3.3 of \cite{WilBd}.  Thus for $A$ in the corresponding closed region, 
$G_A$ is negative on all of $P^\circ$ and so $P^\circ \cap \{ G_A >0 \}$ is empty.  For $E$ a positive multiple of $\tilde B_1$ or $\tilde B_2$, we also have that 
$P^\circ \cap \{ G_E >0 \}$ is empty.  We shall refer to the above cases as the \it trivial \rm cases, as for such classes the results in this section will be satisfied trivially.  

The case of $E$ representing $ B_3$ is more interesting; the conic $\tilde B_3 \cdot D^2 =0$ consists of a line pair with singularity at $R$, one line of which is tangent to $C_2$  there
 (i.e. the asymptote $x+y = k/(k-1)$) and one line of which is the line of symmetry $y=x$.  Thus the subset $Q$ of $P^\circ$ given by $\tilde B_3 \cdot D^2 >0$ lies on one side of the plane through the origin in $\R^3$ 
 corresponding to the line of symmetry, and the arc in $C_2$ given by the same inequality is just the upper half of $C_2$, which we note is also the half which is visible 
  from $\tilde B_3$ with respect to $-Q$.
    Clearly there is a corresponding statement for the point $-\tilde B_3$.  The main cases of interest however are  where $A$ is an affine point, not in the above two closed 
  regions which yield trivial cases, and 
  where the index is $(1,q)$ for some $q\le 2$; here $A$ will be on the tangent line to the Hessian at some point $U$ (not necessarily unique) on  one of the open arcs $Q_1B_3$ or $B_3Q_2$.  The subcone of $P^\circ$ defined by $A\cdot D^2 >0$ will consist of one or two (convex) open subcones, 
  and we wish to describe the (convex) components of $V_2 \cap \{G_A >0 \}$ and  $P^\circ \cap \{G_A >0 \}$.
  
   Given a point $U_1$ on say the open  arc $Q_1B_3$ of the Hessian, we now describe in the Summary below the components of 
 $P^\circ \cap \{ G_A >0 \}$ for classes $A$ on the affine tangent line  at $U_1$ for which the index is $(1,q)$ for some $q\le 2$.  These may also be described 
 in terms of the components of  $V_2 \cap \{ G_A >0 \}$ and $V_1 \cap \{ G_A >0 \}$.
 
\vspace{0.4cm}

\noindent \bf Summary: \rm 
Suppose that $U_1$ is a point of the arc $Q_1B_3$ of the Hessian, with $U_1' = \alpha (U_1)$ the corresponding point of the 
arc $B_1R$ of $C_2$.  The tangent line $L_1$ to the Hessian at $U_1$ will intersect the Hessian again 
at some point $Z$  on the branch 
$B_3B_1$.  Recall that for any point $A$ of the tangent line $L_1$, 
the conic $A\cdot D^2 =0$ contains $U_1'$ and that if $A\ne U_1$, 
then the conic is smooth at $U'_1$ with tangent line $L$ there  being the line joining $U_1$ and $U_1'$. 
 Moreover, 
for any point $A =(a,b,1)$ on the tangent line $L_1$ above the third point of intersection $Z$ with the Hessian, the index at $A$ is $(1,q)$ with $q\le 2$. 
 It follows 
from Proposition 3.4 of \cite{WilBd} that for such points the components of $V_2 \cap \{A\cdot D^2 >0 \}$ 
and $P^\circ \cap \{A\cdot D^2 >0 \}$ are convex.  
We call a component of $V_2 \cap \{A\cdot D^2 >0 \}$ \it bounded 
 \rm if  its boundary intersects $C_2$ in a bounded arc, and \it unbounded \rm otherwise (in which case  its boundary 
 contains $B_1$ or $B_2$).
 
 We consider first the case of $A\in L_1$ having  $a< -{1\over {k-1}}$ 
and $a+b > {k'\over {k'-1}}$ (and so furthest away from the above third point of intersection).  Here $A$ does not lie on any 
other tangent line at a point on the two arcs $Q_1B_3$ and $B_3Q_2$, and so 
$U_1' = \alpha (U_1)$ is the unique point of $C_2$ at which $G_A$ vanishes.  We know that $G_A$ is positive at $B_1$ and vanishes also at  the unique point of $C_1$ for which the tangent passes through $A$.  
 We deduce that $\{ G_A >0\}$ defines a unique component in both $V_2$ and $P^\circ$, the latter containing $-\tilde B_1 = (0,-1,0)$ in its boundary.  The unique component
of $V_2 \cap \{A \cdot D^2 >0 \}$ is an unbounded component and lies above the common tangent line $L$ 
to the conic part of the boundary at $U_1'$, with the corresponding region $V_1 \cap \{A \cdot D^2 >0 \}$ lying below $L$.  
(For the sake of brevity, we shall in future describe components in $V_2$ and leave the reader to formulate the appropriate 
statements about  $V_1 \cap \{A \cdot D^2 >0 \}$.)  Moving down the tangent line $L_1$, when we reach the line $a+b = {k'\over {k'-1}}$, we have that $G_A$ also vanishes (twice) at $R$.

For $A$ between here and $U_1$ on the tangent line, the point $A$ lies on two further tangents, one at $U_2$ in the arc 
$U_1B_3$ and one at $U_3$ on the arc $B_3Q_2$.  With $U_i'$ denoting the corresponding points on $C_2$, we have that $R$ lies 
in the arc $U_2'U_3'$, and that $G_A$ vanishes at these three points on $C_2$, in addition to the point on $C_1$ whose tangent line passes through $A$.  The open 
subset $V_2 \cap \{A \cdot D^2 >0 \}$ then has two components, an unbounded one (with boundary passing through $U_1'$ with $L$ tangent there)
lying above the line $L$,  and a bounded one lying below $L$ whose boundary has points in $C_2$ comprising  the arc $U_2' U_3'$.  

When $A=U_1$, the conic is a real line pair with 
singularity at $U_1'$,
one line passing through the point $\beta (U_1)$ on the arc of $C_1$ between  $B_1$ and the halfway point, determined by the tangent to $C_1$ at $\beta (U_1)$
 passing through $U_1$, and the other line passing through the point $\gamma (U_1)$
on the arc $RB_2$ of $C_2$ which is the image under the Steinerian 
involution of the point on the arc $B_3 Q_2$ for which the tangent to the Hessian passes through $U_1$.  
In particular, we note that the two components of $P^\circ \cap \{U_1 \cdot D^2 >0 \}$ 
 lie on opposite sides of the plane in $\R^3$ corresponding to 
 the line $L$ (for the corresponding points in $V_2$, this means `above' and `below'  the line $L$), since for all $A\ne U_1$ near $U_1$ on the tangent line $L_1$, the quadratic $A\cdot D^2$ is negative on $L\cap V_2$, and thus by continuity the quadratic $U_1 \cdot D^2$ is non-positive on $L\cap V_2$.

Between $U_1$ and the point where $a= -1/(k-1)$, we 
again get an unbounded and a bounded component, but this time it is the bounded component which passes through $U_1'$ 
with tangent $L$ there;
 the unbounded component in $V_2$ (which in this case passes through $U_2'$)  is still above $L$ and the bounded component below $L$.

With $A$ moving further down $L_1$, we have $a\ge -{1\over {k-1}}$, $b \ge -{1\over {k-1}}$ 
and $a+b < {k'\over {k'-1}}$ and we lose the unbounded component.  Here $A$ lies on two tangents to the Hessian, one at $U_1$ on the arc $Q_1B_3$ and one at a point $U_3$ on the arc $B_3 Q_2$,  with images under the Steinerian involution being $U_1'$, $U_3'$, where $R$ is in the arc $U_1'U_3'$ of $C_2$.  In this case the 
inequality $G_A >0$ defines a unique component in $V_2$, lying below the line $L$ and whose boundary has intersection $U_1'U_3'$ with $C_2$, the conic part of the boundary being tangent to $L$ at $U_1'$; 
there is a corresponding unique component in $P^\circ$.   The case $A=Q_2$ is described in more detail below.

As we move further down $L_1$  to points $A$ with $b< -1/(k-1)$, one of two things can happen (until we reach the third point of intersection $Z$ 
with the Hessian --- note here that $Z\ne Q_2$ in this case).  If the $x$-coordinate at $Z$ is less than  the $x$-coordinate at $Q_2$, then $A$ is on a unique tangent line to points on the arc $Q_1B_3$ and $B_3Q_2$, 
namely $L_1$, and so $U_1'$ is the unique point of $C_2$ at which $G_A$ vanishes, with $G_A$ also vanishing at a unique point  
of $C_1$.  Noting that $G_A (B_2) >0$, we have a unique (unbounded) component of $V_2 \cap \{G_A >0\}$ lying below $L$, 
with boundary $G_A=0$ tangent to $L$ at $U_1'$ and with 
 corresponding arc $C_2 \cap \{G_A >0\}$ containing $R$, and hence a unique 
component of $P^\circ \cap \{G_A >0\}$, whose boundary contains $-\tilde B_2$.
 If however the $x$-coordinate at $Z$ is greater than the the $x$-coordinate at $Q_2$, then $A$ lies on the tangents to two 
 points $U_2$ and $U_3$ on the arc $B_3Q_2$, with corresponding points $U_2'$ and $U_3'$ on the arc $RB_2$ of $C_2$.  Here $V_2 \cap \{ G_A >0\}$ has two components, one bounded and one unbounded, with corresponding arcs $U_1'U_2'$ (containing $R$) and $U_3'B_2$ in $C_2$.  Both components lie below the line $L$, and the bounded component has boundary tangent to $L$ at $U_1'$.
 
In the first of the cases above, the conic $Z\cdot D^2 =0$ is a  line pair, which intersects the projectivised boundary of $P$ in 
precisely two points, one of which corresponds to $U_1' \in C_2$ and one to the relevant point of $C_1$ whose tangent line contains $Z$; the conic is smooth at $U_1'$ with tangent line $L$.  Thus one line of the pair must be $L$ and joins $U_1'$ to the relevant point of $C_1$, and the other line is disjoint from $P$.  Therefore $V_2 \cap \{ Z\cdot D^2 >0 \}$ corresponds to the points in $V_2$ below the line $L$.  With $L$ we associate a plane through the origin in $\R^3$, and then $P^\circ \cap \{ Z\cdot D^2 >0 \}$ is the intersection of $P^\circ$ with the 
associated open half-space containing $-\tilde B_2$.    In the case when $Z=Q_2$, the conic $Z\cdot D^2 =0$ is a line pair, with singularity at $B_2$, one line of which only intersects the projectised boundary of $P$ at $B_2$
  and the other line of which is $L$ joining $U_1'$ to $B_2$.  
 Thus again $V_2 \cap \{ Z\cdot D^2 >0 \}$ corresponds to the points in $V_2$ below the line $L$.  Again we associate with $L$ the corresponding plane through the origin in $\R^3$ and $P^\circ \cap \{ Z\cdot D^2 >0 \}$ is the intersection of $P^\circ$ with the open half-space containing $-\tilde B_2$.
 
In the second of the cases above, the conic $Z\cdot D^2 =0$ is a  line pair with singularity at $Z' = \alpha (Z)$ in the open arc 
$RB_2$ of $C_2$.  The conic is smooth at $U_1'$ with tangent $L$ there, so the two lines are $L$ which joins $U_1'$ to $Z'$, and a second line joining $Z'$ to the point on $C_1$ whose tangent contains $Z$.  Here 
$V_2 \cap \{ Z\cdot D^2 >0 \}$ has two components, one being bounded and the other unbounded containing $B_2$ in its boundary.  

 \vspace{0.4cm}

 We now have listed, for the various classes $E$ in the open upper half-space or representing an inflexion point,  with 
 index $(1,q)$ for $q\le 2$, 
 the connected components of $P^\circ \cap \{G_E >0 \}$.  For such a component $Q$,
 we first check that 
  we have the required property that
 the corresponding component $\Gamma$ of  $C_2 \cap \{G_E >0 \}$ is visible from $E$ with respect to $-Q$.
  The crucial general  result needed here is that if the Hessian of an elliptic curve is smooth and $U$ 
 is a point on the Hessian curve  
 with image $U'= \alpha (U)$ under the Steinerian involution, and $U''$ is the third 
 point of intersection of the line $UU'$ with the Hessian, then the tangent lines to the Hessian at $U$ and $U'$ intersect at the point $\alpha (U'') $ of the Hessian (and the line $UU'$ is one of the lines of the line pair 
 $\alpha (U'') \cdot D^2 =0$, with $U''$ being the singularity) --- see \cite{Dolg}, Proposition 3.2.7.  The intersection point 
 of the two tangents is therefore just the third point of intersection of the tangent line to the Hessian at $U$ (or $U'$) with the Hessian.

  \begin{prop}  Given a component $Q$ of $P^\circ \cap \{ E\cdot D^2 >0 \}$, 
  associated to a real  class $E$  in the  extended upper half-space  with 
 index $(1,q)$ for $q\le 2$,
   then all points of $C_2$ whose negative multiples are on the boundary of $Q$ are visible 
  (with respect to the cone $-Q$) from $E$.
  
  \begin{proof}  Recall that by definition $E$ is either in the open upper half-space or represents an inflexion point of the elliptic curve.
  
  This result has already been checked explicitly for the so-called trivial cases, and the cases $E=\pm \tilde B_3$, in the paragraphs preceeding the above Summary.
 We now need to  check visibility for the cases when some positive multiple of $E$ is a point $A$ of the affine 
 plane, not lying in the two closed regions where the result is trivial.  We recall that this then reduced us to considering points 
 $A$, lying on the tangent line to the Hessian at $U_1$, for some $U_1$ on the open arc $Q_1B_3$ or $B_3Q_2$.
 Without loss of generality, we assume that $U_1$ lies on the open arc $Q_1B_3$, and we use the explicit description of the 
 possible components $Q$ from the above Summary.  The quoted classical result says that the tangent $L_1$ at $U_1$ meets the tangent at $U_1'$ at the point $Z$ where $L_1$ meets the Hessian for the third time.  Thus $U_1'$ is visible from $A$ with respect to $V_2$ for all points on $L_1$ above (and to the left of) Z, and is not visible from $A$ for all points on $L_1$ strictly below (and to the right of) Z. 
 
 Let us consider the various possibilities for $A\in L_1$; the first case is $A = (a,b,1) \in L_1$ having  $a< -{1\over {k-1}}$ 
and $a+b \ge {k'\over {k'-1}}$.  Not only is $U_1'$ visible from $A$ with respect to $V_2$, but so also is $B_1$, and thus so 
too is the arc $B_1U_1'$ as required.  The next case will be $a\le -{1\over {k-1}}$ 
and $a+b < {k'\over {k'-1}}$. Here $A$ lies on the tangent line $L_3$ at a point $U_3$ on the arc $B_3Q_2$ of the Hessian  
and $L_3$ 
intersects the Hessian again at a point on the arc $Q_1B_3$ (on the other side of $A$ to $U_3$ on $L_3$). 
 From the above quoted classical result 
 it is this point where the tangent lines at $U_3$ and $U_3'$ intersect.  
 This ensures that $U_3'$ is visible from $A$ as claimed, 
 and thus the same is true for \it all \rm points of the arc $B_1 U_3'$ (not just the points of $C_2$ where $G_A \ge 0$). 
 A similar argument shows that when $a >-{1\over {k-1}}$ and $b> -{1\over {k-1}}$, 
  and so $G_A \ge 0$ defines just a single arc $U_1'U_3'$ in $C_2$, then both endpoints are 
 visible from $A$, as is the arc inbetween.
 
 When $A$ has $b\le -{1\over {k-1}}$ but above the third intersection point $Z$ with the Hessian, we have that $B_2$ is also visible from $A$ with respect to $V_2$, 
 and so too are all points of the arc $U_1'B_2$.
 \end{proof}
\end{prop}

Crucial for the next proof will be the fact noted before, when we introduced the Steinerian map at the start of this Section,  that for any $U$ on the arc $Q_1B_3$ of the Hessian, and $W$ on the tangent line through $U$, not only does the conic $W\cdot D^2 =0$ always intersect $C_2$ at $U' = \alpha (U)$, but also 
when $W\ne U$,
the conic is non-singular at $U'$ and the tangent line to the conic at $U'$ does not depend on the choice of $W$, and is in fact just the line joining $U$ and $U'$.
 The other ingredient that we shall need 
in the proof below is the explicit description of the line pair when the conic is singular.  The Proposition is in fact false if we allow one or both of the $E_i$ to lie in the lower half-space.
 
 \begin{prop} Let $E_i$ (for $i=1,2$) be real classes in the  extended upper half-space with 
 indices  $(1,q_i)$ for $q_i \le 2$, 
  and suppose 
  there are  components $Q(i)$ of $P^\circ \cap \{G_{E_i}>0\}$ for $i=1,2$ whose intersection is non-empty; then some non-trivial  open arc in $C_2$ is in the boundaries of 
 both $-Q(1)$ and $-Q(2)$.
\begin{proof}
Recall that by definition $E_i$ is either in the open upper half-space or represents an inflexion point of the elliptic curve.

If the boundary of one of the components has points corresponding to the whole arc $C_2$, then the result is clear.  So we assume that neither of the $Q(i)$ corresponds to a trivial case.
We now prove that if there is no arc of $C_2$ as described, then $Q(1) \cap Q(2)$ is empty.  Suppose $Q(1)$ is a 
component of $P^\circ \cap \{D\ : \ E_1 \cdot D^2 >0\} $ and $Q(2)$ is a component of $P^\circ \cap
\{ D \  :\  E_2\cdot D^2 >0 \}$.  Corresponding to these components, we have components of $C_2  \cap \{E_1 \cdot D^2 >0 \}$ 
and  $C_2 \cap \{E_2 \cdot D^2 >0\}$, open arcs $\Gamma _1$ and $\Gamma_2$ in $C_2$, where we assume that $\Gamma_1 \cap \Gamma_2$ is empty.
We deduce that at least one component has corresponding arc $\Gamma _i$ in $C_2$ not containing $R$. 
 We comment that if one 
of these arcs has $R$ 
in its closure but not in its interior, then this is the case when $E_i$ is a positive multiple of $\pm \tilde B_3$  and the conic $E_i \cdot D^2 =0$ consists of a line pair with singularity at $R$, one line of which is tangent to $C_2$  there and one line of which is the line of symmetry $y=x$.
Thus the corresponding $Q(i)$ lies on one side of the plane $\Lambda$ in $\R^3$ determined by this line of symmetry. 
The arc corresponding to the other component cannot by our assumption then contain $R$; of course it may be that the other 
component corresponds to taking a negative multiple of $\pm \tilde B_3$, in which case the result is obvious --- although of course we cannot have $\pm \tilde B_3$ \it  both \rm being positive multiples of classes of rigid non-movable surfaces.   Otherwise, 
the argument given in the first basic case below implies that the other component is contained in the complementary half-space, and the result follows.

This enables us to assume that both $E_1$ and $E_2$ lie strictly above the plane $z=0$, and we let $A_1$ and $A_2$ 
denote the corresponding points in the affine plane.
 We can then reduce to considering two  cases: when neither $\Gamma _i$ contains $R$, and when one doesn't and one does.  If $R \not\in \Gamma _i$, we may assume by the above comments that it is not in the closure.
 Using symmetry, the 
above assertion is proved in  these cases from the two basic cases below.  We shall without further reference repeatedly use the facts detailed in the above Summary.

In the first basic case, $Q(1)$ corresponds to an unbounded component of 
$V_2 \cap \{ G_{A _1} > 0 \}$ with the associated arc $\Gamma _1 \subset  C_2$ not containing $R$ and having an endpoint $B_1$ (corresponding 
to a point $A_1$ above the arc $Q_1B_3$ of the Hessian in the region $x < -1/(k-1)$), 
and $Q(2)$ corresponds to an unbounded component of $V_2 \cap \{ G_{A_2} > 0\}$ 
with the associated arc $\Gamma_2 \subset C_2$ not containing $R$ and having an endpoint $B_2$ (corresponding 
to a point  $A_2$ to the right of the arc $B_3Q_2$ in the region $y < -1/(k-1)$).  
We consider the affine picture; if $\bar A_1$ denotes the point of the arc 
$Q_1 B_3$ of the Hessian vertically below $A_1$, we have two lines defined by $\bar A_1 \cdot D^2 = 0$,  
with one line joining $\alpha(\bar A_1)$ to the point $\beta (\bar A_1)$ (\it below \rm the midpoint) 
on $C_1$ where the tangent contains $\bar A_1$, and this corresponds to a plane 
through the origin in $\R^3$.  Recalling that we defined $\tilde B_1 = (0,1,0)$, 
for all $D$ in the corresponding half-space containing $-\tilde B_1$ we 
have $\bar A_1 \cdot D^2 >0$.
Since $\tilde B_1 \cdot D^2 <0$ at all points $D$ of $P^\circ$, the component $Q(1)$ is contained in this half-space. 
A similar statement holds for the component $Q(2)$ ---  we take $\bar A_2$ to be the point on the arc $B_3 Q_2$ of the Hessian horizontally to the left of $A_2$; one  the two lines given by $\bar A_2 \cdot D^2 =0$ (namely 
the one joining $\alpha (\bar A_2)$ to the appropriate point of $C_1$, this point being \it above \rm the midpoint)
corresponds to a half-space in $\R^3$ containing $(-1,0,0) $ and that $Q(2)$ is 
contained in this half-space.  We note that the point of intersection of the two affine lines under consideration is a point of the affine plane not in $V_1 \cup V_2$,
Therefore 
the planes through the origin we have 
constructed via $\bar A_1$ and $\bar A_2$ 
meet in a line disjoint from $P^\circ$, from which it follows that 
$Q(1) \cap Q(2)$ is empty.  In fact, since the above two lines do not intersect the line of symmetry $y=x$ inside $V_1\cup V_2$, 
the two components lie on opposite sides of the plane of symmetry  $\Lambda$ defined above.

The second basic case to consider is when $Q(1)$ corresponds to a component of 
$V_2 \cap \{ G_{A_1}   > 0 \}$, with an associated open arc $\Gamma_1 = B_1 U_1'$ in $C_2$ with $U_1'$ in the open arc $B_1R$ (corresponding to a point $A_1$ 
above the arc $Q_1B_3$ of the Hessian in the region $x < -1/(k-1)$), 
and $Q(2)$ corresponds to a component of 
${V_2 \cap \{ G_{A_2}   > 0 \}}$, with an associated open arc $\Gamma_2$ in $C_2$ containing $R$; 
we saw in the above Summary that $Q(2)$ corresponding to either a bounded component or an unbounded component in $V_2$ was possible here.  We 
 denote the lefthand end of $\Gamma _2$ by $U_2'$, with $U_2'$ either in the open arc $U_1'R$ or equal to $U_1'$.  
The case when $U_2' = U_1' =U'$ is easy, since then the $Q(i)$ lie on opposite sides of the plane in $\R^3$ corresponding to the line joining $U$ to $U'$;
so we assume that $U_2'$ is in the open arc $U_1'R$.  
We have the affine tangent line $L$ to (the conic part of) $-\partial Q(1)$ at $U_1'$, namely the line joining $U_1$ to $U_1'$, and the affine tangent  line $M$ to (the conic part of) $-\partial Q(2)$ at $U_2'$, namely the line joining $U_2$ to $U_2'$.  We note that the point of intersection of the two affine lines under consideration is a point of the affine plane not in $V_1 \cup V_2$, since if  $U$ moves monotonically from $Q_1$ to $B_3$ on the Hessian, its conjugate $U'$ moves monotonically from $B_1$ to $R$ on $C_2$. 
 Moreover, $Q(1)$ corresponds to a subset of the points in $V_2$ above the line $L$ and of the points in $V_1$ below the line $L$, whilst for $Q(2)$ the 
 corresponding points in $V_2$ are below the line $M$ and 
 in $V_1$ are above the line $M$.  By the continuity argument from the Summary, this continues to be true even if one or both of the $A_i=U_i$.
The plane in $\R^3$ corresponding to $L$ determines an open half-space containing $-\tilde B_1$, and we denote by $\tilde Q(1)$ the intersection of this half-space with $P^\circ$, 
and the plane corresponding to $M$ determines an open half-space containing $-\tilde B_2$, and we denote by $\tilde Q(2)$ the intersection of this half-space with $P^\circ$.  
The above planes meet in a line disjoint from $P^\circ$, and so $\tilde Q(1) \cap \tilde Q(2)$ is empty.
Since $Q(i) \subset \tilde Q(i)$ for $i = 1,2$,  the claim follows in this case also.
\end {proof}
\end{prop}

 \begin{cor}  Suppose that the real elliptic curve $F=0$ has two connected  components, then the statement 
 of Proposition 1.10 holds.

  \begin{proof}  For each component $Q(i)$ (with associated class $E_i$) we have an open arc $\Gamma _i \subset C_2$.  
  Suppose first that each $\Gamma_i$ is an arc of the form $B_1 U_i'$; 
   unless some subsequence of the $U_i'$ tends to $B_1$, we have $\bigcap _{i\ge 1} \Gamma _i$ is an arc of the form $B_1U$ in $C_2$ and the result is proved.     Without loss of generality therefore, we may assume that the sequence $U_i'$ tends to $B_1$, and also 
  that no  $E_i$ is a positive multiple of $\tilde B_3$.  Thus we have a point $A_i$ of the affine plane corresponding to each $E_i$. 
 To each $U_i'$, we have a corresponding point $U_i = \alpha (U_i')$ in the arc $B_3Q_1$ of the Hessian, with the $U_i$  tending to $Q_1$, and so the $A_i$ defining $Q(i)$ lies on the tangent line to the Hessian at $U_i$, above (and to the left of) $U_i$.  Moreover for each $i$, we have a common tangent line $M_i$ to the conics 
  $A\cdot D^2 =0$ at $U_i'$ for $A \ne U_i$ on the tangent line at $U_i$, 
  which we saw was just the line joining $U_i$ and $U_i'$.   Corresponding to $M_i$, there is a plane in $\R^3$, and the component $Q(i)$ lies in the associated open half-space that also contains $- \tilde B_1$; as usual, from the continuity argument  in the above Summary,  we note that this is 
  also true when $A_i = U_i$.  We let $\tilde Q(i) \supset Q(i)$ denote the open subcone of $P^\circ$ given by all the points in $P^\circ$ on the same side of the plane corresponding to $M_i$ as $-\tilde B_1$.     If $U_i \to Q_1$, then the tangent  lines $M_i$ tend to the line $M$ corresponding to $Q_1$, which may be checked is 
   the tangent line to $C_1$ at $B_1$, namely the asymptote  given affinely 
  by $x =-1/(k-1)$; it intersects the projectivised boundary of $P$ just at $B_1$.  This implies that 
  $\bigcap _{i\ge 1} \tilde Q(i)$ is empty, and hence so too is $Q = \bigcap _{i\ge 1}  Q(i)$, contrary to assumption.  Similarly, we see that the result is true when each $\Gamma _i$ has $B_2$ as a righthand end.

  We now put an ordering on the points of $C_2$ by specifying that $S \le S'$ if $S$ is in the (closed) arc $B_1 S'$.  Let $S_1 \in C_2$ denote the infimum of the 
 righthand ends 
 of the arcs $\Gamma _i =  \partial (- Q(i)) \cap C_2$
 and let $S_2 \in C_2$ denote the supremum of the lefthand ends of arcs  $\Gamma _i $.
 
 If $S_2 < S_1$, then the required open arc of $C_2$ is $S_2S_1$.
  If $S_2 > S_1$, then we can find components $Q(1)$ and $Q(2)$ say, 
 with $\Gamma _1 \cap \Gamma_2$ empty.  From the Proposition 2.2, we then have that $Q(1) \cap Q(2)$ is empty, therefore contradicting the assumption that $Q$ has non-empty interior. 

   \begin{figure}
 \centering
     \includegraphics[width=8cm]{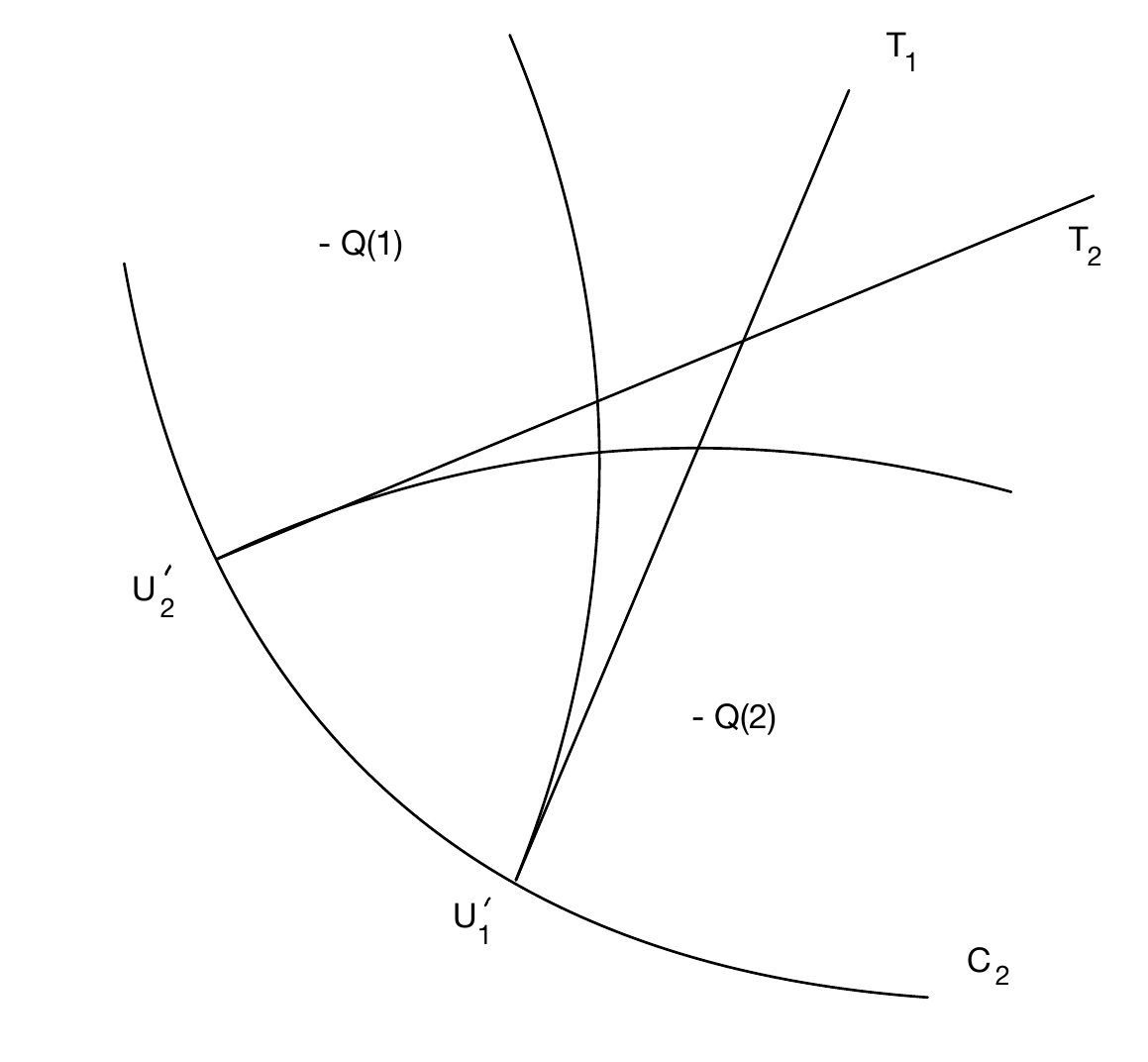}
\caption{Diagram for last part of the proof of Corollary 2.3}
%\label{}
\end{figure}

 Finally, we need to deal with the remaining possibility, that for any pair of components $Q(1)$ and  $Q(2)$ say, including the case of $Q(1) = Q(2)$ corresponding to a bounded component in $V_2$, 
 we have that $\Gamma _1 \cap \Gamma _2$ is an
arc $U_2' U_1'$ of $C_2$ in the boundary of  $-Q(1) \cap Q(2)$, but the intersection of all such arcs is a single 
point $U' = \alpha (U)$ on $C_2$.  We may assume without loss of generality 
that $U_2' \ne B_1$ is a lefthand end of $\Gamma_2$ and $U_1' \ne B_2$ is a righthand end of $\Gamma _1$, and 
we have a schematic diagram as in Figure 2,  
illustrating the regions of $V_2$ corresponding to $-Q(1)$ and $-Q(2)$.  For $i=1, 2$, we saw that for any $A_i$ on the tangent to the Hessian at $U_i$, the conic $A_i \cdot D^2 =0$ has a zero at $U_i' = \alpha (U_i) \in C_2$ and the tangent $T_i$ to the conic at $U_i'$ is independent of the choice of $A_i \ne U_i$, 
and was in fact just the line joining $U_i$ and $U_i '$.  It follows that all the points of $V_2$ which are negative multiples of elements in $Q(1)$, respectively $Q(2)$, lie above $T_1$, respectively below $T_2$, whilst the points of $V_1$ which are positive multiples 
of elements in $Q(1)$, respectively $Q(2)$, lie below $T_1$, respectively above $T_2$.  This ensures that $Q(1)\cap Q(2)$ is
contained in a convex subcone of  $P^\circ$ lying between the planes corresponding to $T_1$ and $T_2$.  By the 
continuity argument from the Summary, this continues to be true even if one or both of the $A_i = U_i$.  In the case when $Q(1) = Q(2)$ is a bounded component, the statement 
also remains true.

If now $L$ is the common tangent line $UU'$ to the conics $A\cdot D^2 =0$, where $A\ne U$ lies on the tangent to the Hessian at 
the point $U = \alpha (U')$ defined above, then we can find a subsequence of  the $Q(i)$ for which the corresponding tangent lines $T_1$ at the righthand ends of $\Gamma _i$ tend to $L$, and similarly a subsequence of the $Q(i)$ for which the corresponding tangent lines $T_2$ at the lefthand ends also tend to $L$.  From this we deduce that the points of $Q = \bigcap _{i\ge 1} Q(i)$  must lie 
in the plane through the origin  in $\R^3$ corresponding to $L$, and so $Q$ would have empty interior, contrary to assumption. 
 \end{proof}
\end{cor}

Therefore in the case when the real elliptic curve has two components, we have completed (via Proposition 1.9) the proof of Theorem 1.8,  and hence 
we have proved (via Corollary 1.7) the relevant parts of our
Theorem 0.2.

\section{Hybrid components when elliptic curve has one real component}

 We now wish  to study the case when the elliptic curve has only one real component, and so the Hessian,
 assumed smooth, 
has two components. 
 Our objective will be to prove the analogous results to Proposition 2.1, Proposition 2.2 and Corollary 2.3 in this case, and in this way complete (via Proposition 1.9) the proof of Theorem 1.8, and hence prove the relevant part of our Theorem 0.2.

 The three hybrid components
 of the positive index cone are  as 
described in Section 1.  Recall that there are two special values for $k<1$ where changes occur, namely $k =0$ and $-2$.  Away from these two values, we wish to describe the Steinerian map $\alpha$.

  If as before the inflexion points of the 
cubic (and hence of the Hessian) are denoted $B_1, B_2, B_3$, then the tangents there  to 
the cubic $F$ (which we saw are just the asymptotes to the three affine branches) will be 
tangent to the Hessian at three distinct points $Q_1, Q_2, Q_3$.  These lines may also be characterised as second polars of $F$ with respect to the $B_i$. 
 Having chosen one of the $B_i$ as the zero of the group law, the corresponding point $Q_i$  where the tangent to $F$ at $B_i$ is tangent to the Hessian 
is just one of the real  2-torsion points of the Hessian.   

It is easiest to understand what is going on dynamically.  
For $k>1$, we found a  description of $\alpha$, 
where the tangent to $F$ at each $B_i$ is tangent to an affine branch of the unique connected component of $H$.  If we consider the  corresponding points of the upper half-sphere in $S^2$,  we note that as $k\to 1$, 
 the affine branches of both 
the real curves given by $F$ and $H$ tend to arcs of the equator $z=0$ between representatives of the  relevant inflexion points, 
whilst the bounded component shrinks to a point, so that for $k=1$ both $F$ and $H$ vanish on the equator plus an isolated point on the upper half-sphere 
corresponding to the centroid $({1\over 3}: {1\over 3}:1)$ of the triangle of reference.  Deforming away from $k=1$ towards zero, the singular  point then 
expands to become the 
bounded component of the Hessian, and the arcs on the equator that were limits as $k\to 1+$ 
of the affine branches of $F$ deform to affine branches of $H$ and  the arcs that were limits as $k\to 1+$ 
of the affine branches of $H$ deform to affine branches of $F$.  Recall here from Section 1 that $H_k = -54k^2 F_{k'}$ where $k'  = \frac{4-k^3}{3k^2}$, 
and so one does expect the regions of the affine plane  occupied 
by the unbounded affine  branches of $F$ and $H$ to switch over.
 In particular, for $0<k<1$, the tangents to the cubic at each inflexion point are tangents at appropriate points 
 to the unbounded affine branches of $H$, similar therefore to the case $k>1$.  
Thus in this case, the Steinerian map $\alpha$ sends each inflexion point $B_i$ to a point on the unbounded component of $H$, and hence gives an involution on both 
connected components of $H$ individually.  

The next change occurs at $k=0$, where the bounded component together with the three unbounded affine branches of $H$ just tend to the three
real lines determined by the triangle of reference.  To see what happens to the Steinerian map, it is probably easiest to look instead 
at the value $k= -2$; here the Hessian is just the line at infinity together with the isolated point $({1\over 3}: 
{1\over 3} :1)$ and all three asymptotes of the cubic pass through this point.  As one deforms in either direction away from $k=-2$, this point expands to give the bounded component of the Hessian and each asymptote of $F$ will 
now be tangent to the bounded component of the Hessian, which by continuity will also be the case for all 
$k<-2$ and $-2 < k <0$.  Thus for $k < -2$ and $-2 < k < 0$, the Steinerian map $\alpha$ interchanges the two components of $H$.  The bounded component will be contained in (and tangent to) the asymptotic triangle given by the lines $x ={1\over {1-k}}$, $y = {1\over {1-k}}$ and $x + y = {{k}\over {k-1}}$.  For $k >-2$, the asymptotic  triangle will be given by the inequalities $x \le{1\over {1-k}}$, $y \le {1\over {1-k}}$ and $x + y \ge {{k}\over {k-1}}$, whilst for $k<-2$, it will be given by $x \ge {1\over {1-k}}$, $y \ge {1\over {1-k}}$ and $x + y \le {{k}\over {k-1}}$.

What is occurring here is that for each value of $k' > 1$, there are three possible values of $k$ for which $H_k$ is a multiple of $F_{k'}$, one with $k<-2$, one with $-2 < k < 0$ and one with $0< k <1$.  If we choose an inflexion point  $B_3$ say, the tangent to $F_k$ at $B_3$ will be tangent at 
one of the three real 2-torsion points of $H_k$ and hence $F_{k'}$, the one on the unbounded component if $0< k <1$, and the ones on the bounded component in the other two cases.   

For any given point $A= (a,b,1)$ in the affine plane $z=1$ at which the index is $(1,q)$ with $q\le 2$, we let $G_A$ denote the homogeneous quadratic given by $A\cdot D^2$, explicitly 
in coordinates $$-ax^2 -b y^2 - (1-a-b)(z-x-y)^2 +kay(z-x-y) +kbx(z-x-y) + k(1-a-b)xy, $$
and we wish to understand how 
$G_A =0$ intersects not only the (unbounded)  affine branches of $F$ but also the unbounded affine branches of the Hessian.   We will therefore need to understand this in all the three cases detailed above, as the Steinerian map will be different in the three cases.

In all three cases, we let 
 $C_1$ denote the unbounded branch of $F=0$ which lies in the region 
$x>0, \ y>0$ and $x+y >1$, and $C_2$ the unbounded branch of the Hessian lying in the sector $x<0,\  y<0$.  
There is then a hybrid component $P^\circ$ of the positive index cone whose boundary consists of the positive cone on $C_1$ together with the negative cone on $C_2$, the two parts meeting along rays generated by $(0,1,0)$ and $(1,0,0)$, and without loss of generality we may assume that this is the hybrid component which we study.

 \begin{figure}
 \centering
     \includegraphics[width=10cm]{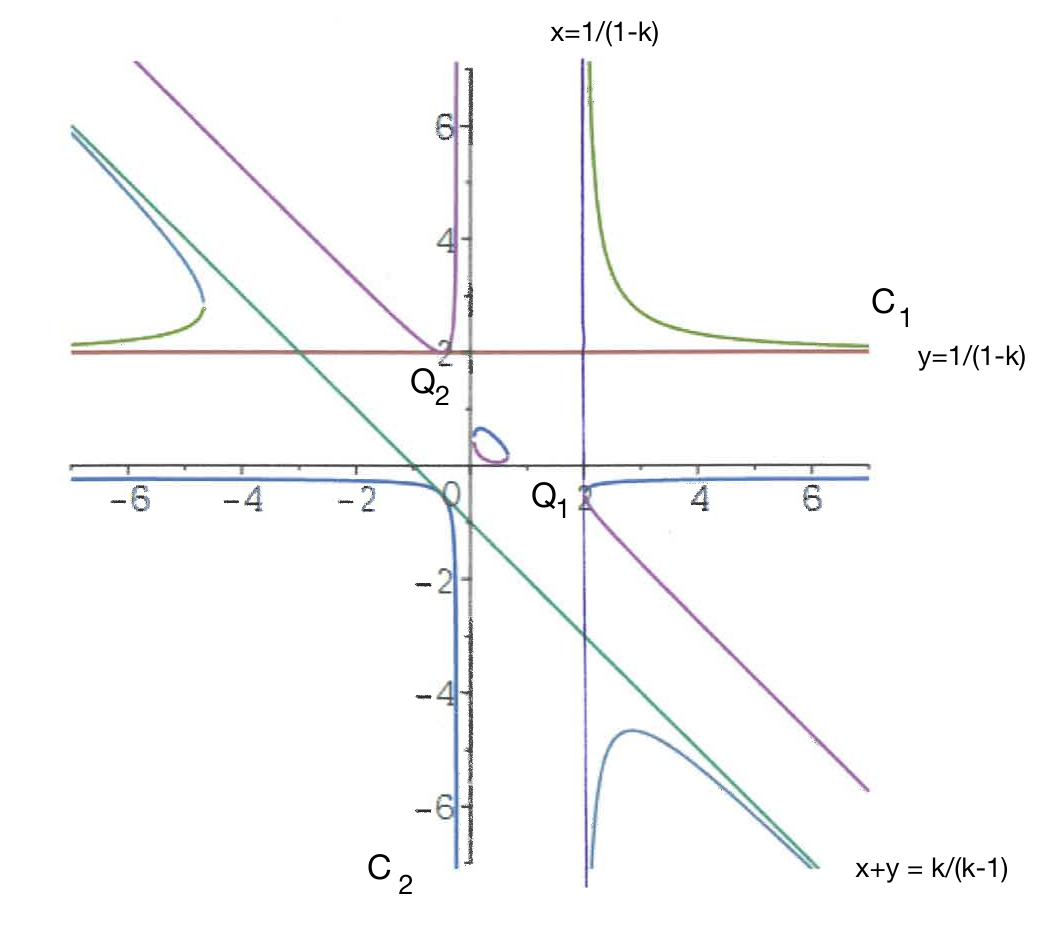}
\caption{Cubic $F_k =0$ with asymptotes and Hessian $H_k$ for $k=0.5$}
%\label{}
\end{figure}

For the case when the cubic is $F_k$ with $0<k<1$, the proofs of the analogous results to those in Section 2 
 are essentially identical to the arguments for $k>1$ given in Section 2, modulo the fact that the 
regions of the affine plane occupied respectively  by the unbounded branches of $F$ and $H$ have switched over.  This case is illustrated in Figure 3 (which shows $F_k$, $H_k$ and the three asymptotes for $F_k$ when $k= 0.5$).  Here the Hessian consists of a bounded component, $C_2$  and the components containing $Q_1, Q_2$, whilst the real cubic curve consists of the other three unbounded components (including $C_1$) as illustrated.
The tangent to $F_k$ at $B_i$ (i.e an asymptote) is therefore tangent to the unbounded component of the Hessian at $Q_i$, for $i=1,2,3$.  As in the case of two components, the asymptote $x+y = k'/(k'-1)$ to the Hessian will again play a pivotol role.  If an affine point $A$ is in the closed region $x\le 1/(1-k)$, 
$y\le 1/(1-k)$ and $x+y \le k'/(k'-1)$, then $G_A$ is negative on $P^\circ$.  If $A$ is in the closed quadrant $x\ge 1/(1-k)$, $y\ge 1/(1-k)$, then $G_A$ is positive on $C_2$ and all of $C_2$ is visible from $A$ with respect to $V_2$. At all other affine points with index $(1,q)$ for some $q\le 2$, 
 we have that $A$ is on the tangent line to the Hessian at some point $U$ (not necessarily unique) on one of the open arcs $Q_2B_3$ or $B_3Q_1$.  Under the Steinerian involution, the arc $Q_2B_3$ corresponds to the arc $B_2Q_3$ of $C_2$ and the arc $B_3Q_1$ to the arc $Q_3B_1$ of $C_2$.
Analogous  arguments to those in Section 2 may then be seen to prove Proposition 1.8 in the case $0<k<1$.   For $k>1$, 
the bounded component of $F$ essentially played no role in the proof of Proposition 1.10, whilst for $0 < k <1$,  it is 
the bounded component of $H$ that essentially plays no role in the proof.  We shall therefore 
not give any further details  in this case, and from now on concentrate on the other two cases.  Even in these other two cases, although the Steinerian involution looks very different,  the arguments we use are very similar to those in Section 2.

\vspace{0.5cm}
\it For the rest of this Section, we shall assume that $-2 < k <0$, and in the next Section we shall study the remaining case when $k<-2$.  We shall in this Section prove Proposition 1.8 when $-2<k<0$.
\vspace{0.5cm}
\rm

For $-2 < k <0$, we have a (schematic) picture as in Figure 4, showing all the affine branches of the Hessian and the branch $C_1$ of the cubic. 
We note that the tangent line to $F$ at $B_3$ is the line 
$x+y = e_2 = {k\over {k-1}} = {{-k}\over {1-k}}$, and this is tangent to  the bounded component of the Hessian at the point $Q_3 = (e_2/2 , e_2/2)$.  
If we take $B_3$ as the zero of the group law, $Q_3$ is just a 2-torsion point of $H$.  Moreover the bounded component of the Hessian
  in this case lies above (and touches) the asymptote ${x+y }= {k\over {k-1}}$.  
Also playing a role will be the other two lines through $B_3$ that are tangent to the Hessian and yielding 2-torsion points of the Hessian; these have the form $x+y = e_3$, 
corresponding to the other tangent to the bounded component  and $x+y = e_1$ corresponding to the tangent to the unbounded component of $H$.  
Explicitly, if the Hessian is  (up to a multiple) the Hessian of $F_{k_i}$, 
where $0 < k_1 <1$, $-2 < k_2 = k  <0$ and $k_3 < -2$, then $e_i = k_i/(k_i -1)$; moreover $e_1 < 0  < e_2 < e_3 < k'/(k'-1)$, where $x+y = k'/(k'-1)$ is the asymptote to the 
Hessian corresponding to its tangent line at $B_3$.

  \begin{figure}
 \centering
     \includegraphics[width=10cm]{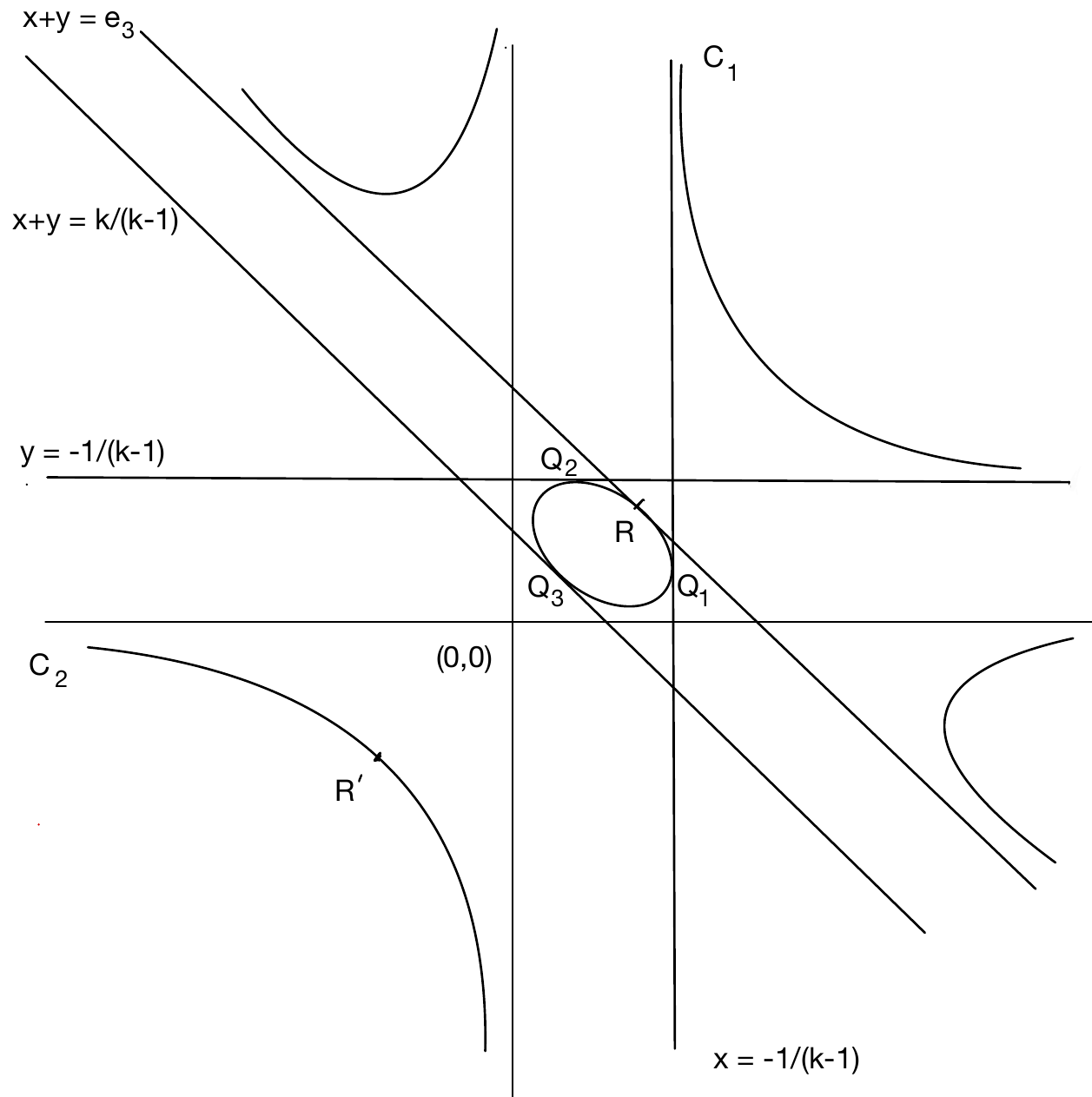}
\caption{Schematic picture when $-2 <k<0$}
%\label{}
\end{figure}

We have the Steinerian involution on the Hessian which in the case being studied interchanges the two components; explicitly we 
let $Q_1$ be the point on the bounded component of the Hessian whose tangent is also the tangent to $F$ 
at $B_1$, and $Q_2, Q_3$ defined similarly with respect to the inflexion points $B_2, B_3$.  We saw 
above that $Q_3 =(e_2/2 , e_2/2)$ and the tangent line is $x+y =e_2$, where $e_2 = {k\over {k-1}}$. 

Under the Steinerian map, the branch $C_2$ (going from $B_1$ to $B_2$)  
of the Hessian corresponds to the `upper'  arc (i.e. not containing $Q_3$) on the bounded component going from $Q_1$ to $Q_2$.  We now argue similarly to Section 2.  For a given class $A = (a,b,1)$,  it is clear how many times the conic $G_A = A\cdot D^2 =0$ cuts $C_1$ --- if $a\ge {1\over {1-k}}$ and $b\ge {1\over {1-k}}$, 
then either $A\in V_1$ and so $G_A>0$ on both $C_1$ and $P^\circ$, or $G_A =0$ intersects $C_1$ twice
 (since there will be two tangents to $C_1$ containing $A$, including maybe at points at infinity or the degenerate case $A\in C_1$);  it will not meet $C_1$ at 
all if $a< {1\over {1-k}}$ and $b< {1\over {1-k}}$, and will cut it precisely once in the other cases.  
We now ask how many times and where the conic cuts $C_2$.  To answer this question, we are looking for the tangents from $A$ to 
the upper arc (i.e. not containing $Q_3$)  from $Q_1$ to $Q_2$ 
on the bounded component of the Hessian.  Here the answer is twice (counted with multiplicity) if 
$A$ is in the region bounded by $a 
= {1\over {1-k}}$, $b ={1\over {1-k}}$ and by the specified arc  $Q_1 Q_2$, none for any  other points with 
$a< {1\over {1-k}}$ and $b< {1\over {1-k}}$ or with $a> {1\over {1-k}}$ and $b> {1\over {1-k}}$, and 
precisely once otherwise as there is exactly one tangent to the given arc $Q_1 Q_2$.  If  $a\ge {1\over {1-k}}$ and $b\ge {1\over {1-k}}$, then $G_A$ is positive on the affine branch $C_2$, all of which is visible from $A$.  If  $a\le {1\over {1-k}}$ and $b\le {1\over {1-k}}$, and $A$ is not in the region specified above, then the same continuity argument we used in Section 2 shows that $G_A$ is strictly negative in $P^\circ$ since it is when $A\in V_2$.

The midpoint of 
the arc $Q_1Q_2$ is the point $R = (e_3/2 , e_3/2)$, and under the Steinerian map this point corresponds to the midpoint 
$R' = \alpha (R) = (e_1/2 , e_1/2)$ 
of $C_2$, namely the intersection of $C_2$ with $x=y$.  Thus $G_A =0$ will cut $C_2$ in the part given by $y<x$ if and only if $A$ lies on a tangent to the Hessian at some point on the open arc $RQ_1$.

The reader will easily check that the points $E= \pm \tilde B_i$ for $i=1,2$ give rise to trivial cases where either $G_E >0$ on both $C_2$ and $P^\circ$ and all of $C_2$ is visible from $E$, or $G_E <0$ on $P^\circ$.
One difference from the first two cases considered for the elliptic curve is that there, the points $E = \pm \tilde B_3$ were isolated, in that they did not represent the point at infinity on any tangent to the affine Hessian curve.  In the remaining two cases, this is no longer true, in that they represent the point at infinity of the tangent line at $R$, the midpoint of the 
relevant arc $Q_2Q_1$.  It is still the case that 
$E\cdot D^2 =0$ is a line pair, with singularity at $R' = \alpha (R) \in C_2$, one line being tangent there and the other being the line 
of symmetry; in the two cases now being studied,  the points $A$ on the affine tangent line give rise to conics $A\cdot D^2 =0$ containing $R'$, and that for $A\ne R$ the conic is smooth at $R'$ with tangent line $L$  the line of symmetry.  The special 
cases $E = \pm \tilde B_3$  are not then as special as previously, in that now they represent a limit of points on the affine tangent line at $R$.  The proofs can therefore invoke continuity when dealing with these cases.

Let us consider a point $U$ on the open arc $Q_2R$; the tangent at $U$ will intersect the Hessian again on 
the branch $B_2B_3$.  We assume that $A$ lies on this tangent line and the index at $A$ is $(1,q)$ with $q\le 2$.  Thus $A$ lies above (and to the left of) the 
third intersection point
$Z$ with the Hessian, and for any such $A$, the conic $A\cdot D^2 =0$ passes though the point $U' = \alpha (U)$ on the arc $B_2R'$ of $C_2$.  We have in this case that $Q = P^\circ \cap \{ G_A >0 \}$ is always connected.  For the point $Z$, the line pair $Z\cdot D^2 =0$ intersects the projectivised boundary 
of $P$ in two points, 
one on $C_2$ (namely $U'$) and one on $C_1$.  Moreover it is non-singular at $U'$ with tangent line $L$, and so one of the two lines must be $L$ and the other will be disjoint from the projectivised boundary of $P$.  The line $L$ may then be identified not only as the line joining $U$ to $U'$, but also as the line joining $U$ to the point 
$U''$ on $C_1$ whose tangent line  passes through $Z$.
  As $U$ moves from $Q_2$ to $R$ monotonically, $U'$ moves monotonically  from $B_2$ to $R'$, Z moves monotonically from $B_2$ to $B_3$ and $U''$ moves from $B_2$ to the midpoint of $C_1$. This latter progression is not however monotonic if $F_k (  \frac{k'}{2(k'-1)}, \frac{k'}{2(k'-1)}, 1) <0$ (this is the case illustrated schematically in Figure 4 and  occurs for example when $-\frac{1}{2} \le k <0$) --- under the given  condition, the tangent to $C_1$ at its midpoint intersects the branch $B_2B_3$ of the Hessian at $B_3$ and a further point.    
 This is the reason why a simple-minded argument just involving tangent lines to the conics does not suffice to prove Proposition 3.2, unlike in our later proof of Proposition 4.2.
 Note that  $ P^\circ \cap \{ G_Z >0 \}$ is given by intersecting $P^\circ$ with an open  half-space corresponding to the plane through the origin  determined by the line $L$. 
 
If $A= (a,b,1)$ above (and to the left of) $Z$
with $a \ge {1\over {1-k}}$, then the closure of $Q$  contains 
$\tilde B_1 = (0,1,0)$, and in the case of strict inequality it contains not only $(0,1,0)$ but also points of $C_1$.
When $a <{1\over {1-k}}$, we have that $(0,1,0)$ is no longer in the closure of $Q$.  When $A$ 
is between the point with $a ={1\over {1-k}}$
and $U$, there is a second point $U_1$ on the arc $UQ_1$ of the Hessian for which $A$ also lies on the tangent line at $U_1$; thus there is a point $U_1'$ in the open arc $ U' B_1$ at which $G_A$ vanishes.  In this case $\partial (-Q) \cap C_2$ consists of the finite arc 
$U' U_1'$.  For $A= U$, the conic is a pair of \it complex \rm lines with singular point $U'$, and $G_A <0$ on $P^\circ$.
As $A = (a,b,1)$ passes to the other side of $U$ but with $b < {1\over {1-k}}$, we again obtain a second zero 
$U_2'$ of $G_A$ on $C_2$, this time in the arc $B_2U'$, corresponding to the other point $U_2$ in the arc $Q_2 U$ 
where the tangent contains $A$, and $\partial (-Q) \cap C_2$ consists of the finite arc $U_2'U'$.  When $b = {1\over {1-k}}$, then the closure of 
$Q$ contains $\tilde B_2 = (1,0,0)$, and then for all further points $A$ we have that the closure of $Q$ contains points of $C_1$ in addition to $(1,0,0)$ 
in its boundary.

 We now use similar methods as in the previous section to prove analogous results to Proposition 2.1, 
 Proposition 2.2 and Corollary 2.3.

  \begin{prop}  Suppose the elliptic curve has one component, with invariant $-2 < k<0$.
  Given a component  $Q$ of $P^\circ \cap \{ E\cdot D^2 >0 \}$, associated to a real class $E$  in the  extended upper half-space  with 
 index $(1,q)$ for $q\le 2$,
  all points of $C_2$ whose negative multiples are on the boundary of $Q$ are visible 
  (with respect to the cone  $-Q$) from $E$.
  
\begin{proof} Since visibility is a closed condition, we may assume that $E$ represents a point 
$A=(a,b,1)$ in the affine plane.  
Moreover we may assume that $A$ does not lie in either of the two regions specified above where the result is trivially true.   
If $A$ lies on the tangent to the Hessian at $R$, strictly above (and to the left of) the point $R$, then it is the upper half $B_2 R'$ of $C_2$ on which $G_A$ is positive, and all these points are visible from $A$ (it is only when $E = \tilde B_3$ that this is precisely the set of 
points visible from $E$).   By symmetry,  the result holds also when $A$ lies on the tangent to the Hessian at $R$, strictly below (and to the right of) the point $R$.

Otherwise, we may assume without loss of generality that $A$ lies on the tangent line at 
a point $U$ in the open arc $Q_2R$, since then the case when $U = Q_2$  follows by continuity.
We noted above, that for such points $A$, the tangent line 
at $U$ intersects the Hessian at a point on the branch $B_2B_3$.    We note that $U'  = \alpha (U)$ lies on the arc $B_2R'$ of $C_2$.

 By the classical result
used repeatedly in Section 2 (\cite{Dolg}, Proposition 3.2.7), the point where the tangents to the Hessian at $U$ and $U'$ meet is precisely the third point of intersection of the tangent at $U$ (or $U'$)  with the Hessian.   Given that the Hessian at $A$ is non-negative, $A$ 
lies above (and to the left of) this point of intersection $Z$, from which it follows that $U'$ is visible from $A$.

As  the point $A= (a,b,1)$ moves up the tangent line from $Z$
 towards $U$, the part of $C_2$ on which $G_A >0$ 
is initially just the arc $U'B_1$, but this is contained in the arc of points which are  visible from $A$. When we reach points $A$ with $a < {1\over {1-k}}$, the part of $C_2$ on which $G_A >0$ is then a
bounded arc $U'U_1'$, whilst all of $C_2$ is visible from $A$.  The case $A=U$ is  not relevant here, and 
as $A$ passes to the other side of $U$, 
the part of $C_2$ on which $G_A >0$ 
is then initially a bounded arc $U_2'U'$, and when $b \ge {1\over {1-k}}$ this is all of $B_2U' $.  With $A$ moving 
upwards from $U$ on the tangent line, 
initially all of $C_2$ is visible from $A$, but the set of visible points will always be 
an arc containing $B_2U'$; hence  all points of the arc $B_2U'$ are visible from $A$, which verifies the claimed result.
\end{proof}
\end{prop}

 \begin{prop}  Suppose the elliptic curve has one component, with invariant $-2 < k<0$.
 Let $E_i$ (for $i=1,2$) be real classes  in the  extended upper half-space with 
 indices $(1,q_i)$ for $q_i\le 2$,
   and suppose 
  there are  components $Q(i)$ of $P^\circ \cap \{G_{E_i}>0\}$ for $i=1,2$ whose intersection is non-empty; then some non-trivial  open arc in $C_2$ is in the boundaries of 
 both $-Q(1)$ and $-Q(2)$.

 \begin{proof} As in the proof of Proposition 2.2, we show that if 
  the arcs on $C_2$ corresponding to the $Q(i)$ are disjoint, then  $Q(1) \cap Q(2)$ is empty.  
   Suppose $Q(1)$ is $P^\circ \cap \{D\ : \ E_1 \cdot D^2 >0\} $ and $Q(2)$ is $P^\circ \cap
\{ D \  :\  E_2\cdot D^2 >0 \}$.  Corresponding to these components, we have open arcs $ \Gamma_1 =C_2  \cap \{E_1 \cdot D^2 >0 \}$ 
and  $\Gamma _2 = C_2 \cap \{E_2 \cdot D^2 >0\}$ in $C_2$. 
The assumption that $\Gamma _1$ and $\Gamma_2$ are disjoint, means that for the appropriate endpoints $U_1'$ and $U_2'$ 
of $\Gamma _1$ and $\Gamma _2$ on $C_2$, we may without loss of generality assume that $\Gamma _1$ is a sub-arc of 
$B_2U_1'$ and that $\Gamma _2$ is a sub-arc of $U_2'B_1$ with $U_1'$ in the arc $B_2U_2'$. 
 We then have corresponding points $U_1$, $U_2$ on the arc 
 $Q_2 Q_1$ of the bounded component ($U_1$ in the arc $Q_2 U_2$).  We shall also assume that the $E_i$ determine points $A_i$ in the affine plane --- as we note below, the limit cases with $E_i $ a positive multiple of $\pm \tilde B_3$ (when the corresponding $U_i =R$) will follow by an essentially unchanged argument.  
 Thus $Q(1)$ corresponds to a point $A_1$ on the tangent line $L_1$ at $U_1$, with $A_1$ strictly above (and to the left of) $U_1$, whilst $Q(2)$ corresponds to a point $A_2$ on the tangent line $L_2$ at $U_2$, with $A_2$ strictly below (and to the right of) $U_2$. With these conventions, we 
 show that $Q(1)\cap Q(2)$ is empty.  Recall that for all $W\in L_1$ the conic $W\cdot D^2 =0$ passes though 
 $U_1' \in C_2$, and if $W \ne U_1$ the conic is smooth there with tangent line $L$ independent of choice of $W$, specifically $L$ is the line joining $U_1$ to $U_1'$.  Moreover $-Q(1) \cap V_2$ lies above $L$ in $V_2$ (and $Q(1)\cap V_1$ lies below $L$ in $V_1$).

\begin{figure}
 \centering
     \includegraphics[width=8cm]{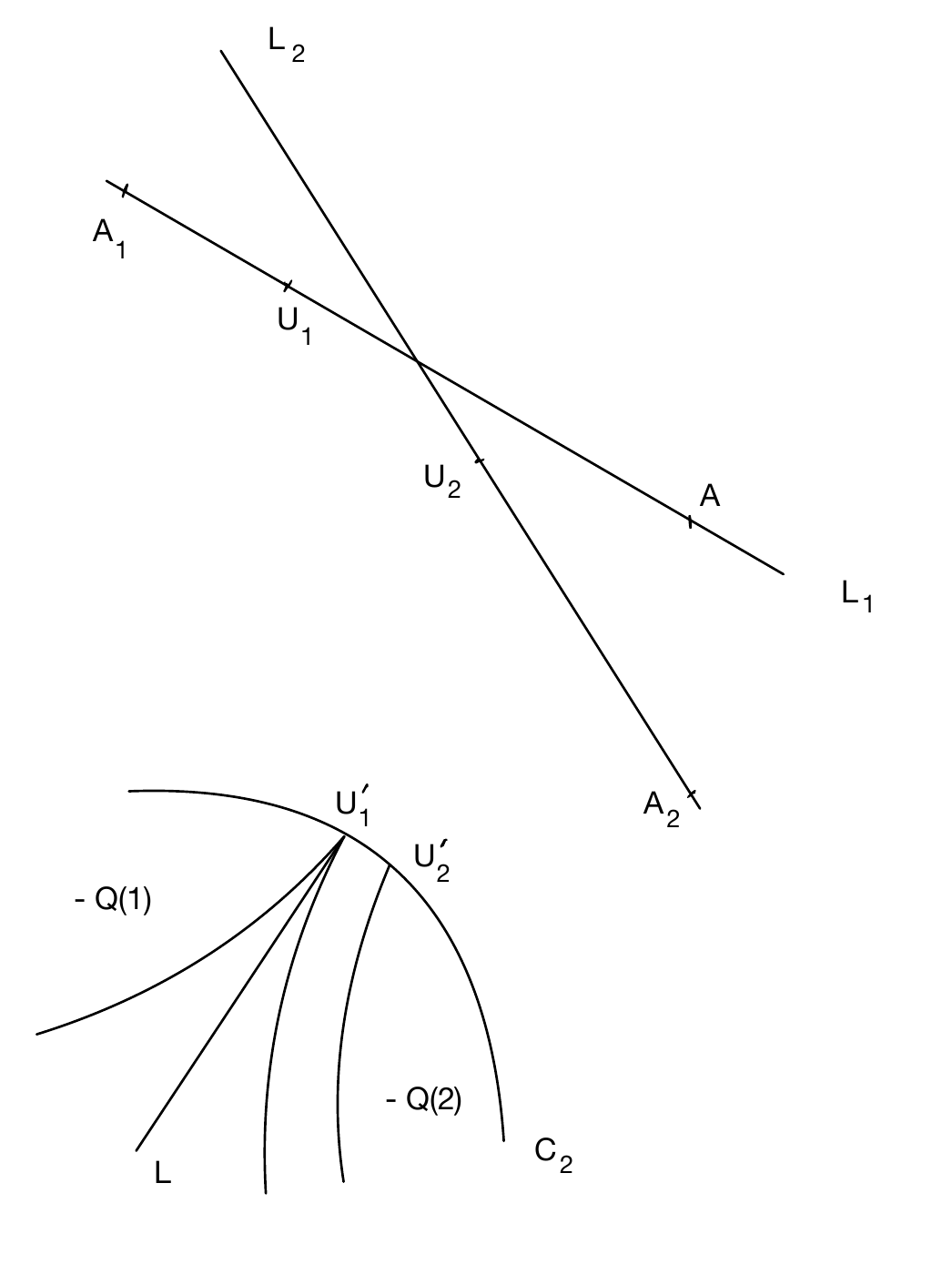}
\caption{Schematic diagram for the proof of Proposition 3.2}
%\label{}
\end{figure}

 We consider various possibilities for $A_1$ and $A_2$; the easy case is when $U_1 = U_2$, for then $-Q(2) \cap V_2$ lies below the 
 common tangent line $L$ in $V_2$ (and $Q(2)\cap V_1$ lies above $L$ in $V_1$), and disjointness of the two components is clear.  
 Let us assume then by symmetry that $U_2$ is in the arc  $RQ_1$.
 We denote by $A$ the point on $L_1$ vertically above $A_2$ as in the main diagram in Figure 5, and hence below (and to the right of)  $U_1$  on $L_1$ (in fact also below the intersection point of $L_1$ and $L_2$).  Suppose first that the Hessian is non-negative at $A$.
   If we consider the subcone of $P^\circ$ given by $A \cdot D^2 >0$, this (projectively) lies on the opposite side of $L$ to that given by $A_1 \cdot D^2 >0$; thus the 
 subcones of $P^\circ$ given by $A_1 \cdot D^2 >0$ and $A \cdot D^2 >0$ lie on opposite sides of the plane in $\R^3$ corresponding to $L$, and hence in particular are disjoint.
 Since $\tilde B_1\cdot D^2 >0$ for all $D \in P^\circ$, we know that the 
open connected subcone of $P^\circ$ corresponding to $A_2$ is strictly smaller than the 
relevant subcone corresponding to $A$, and hence the result follows in this case.  This 
argument is illustrated schematically in Figure 5, where the main diagram shows the configuration of points and tangent lines to the arc $Q_2Q_1$ on the bounded component of the Hessian, where in order not to complicate the diagram we have omitted the arc (which has tangent $L_1$ at $U_1$ and $L_2$ at $U_2$). The smaller subdiagram shows the corresponding regions in $V_2$.  Here $L$ is not only tangent to $-Q(1)$ at $U_1'$, with $-Q(1)\cap V_2$ lying above $L$, but it is also tangent  to the region defined in $V_2$ by $A\cdot D^2 >0$, shown in the diagram lying below $L$ and containing $-Q(2) \cap V_2$.
  In the case where one of the $U_i =R$, without loss of generality $U_1 =R$ and $U_2$ is in the open arc $RQ_1$, the same proof works, even for the limit case when $E_1 = \tilde B_3$, where $Q(1)$ consists of the points of $P^\circ$ lying on the appropriate side of the plane of symmetry.
  
  As a special case, we note that the result  is true when both $U_1$ and $U_2$ lie in the (closed) arc $RQ_1$, since the Hessian at $A$ is 
  then clearly positive.  By symmetry, the result is also true when 
   both $U_1$ and $U_2$ lie in the (closed) arc $Q_2R$; if we write out the analogous proof to the one above, we shall be considering the point given by the horizontal projection of $A_1$ onto $L_2$.

We are left therefore with the case when $U_1$ is in the open arc $Q_2R$ and $U_2$ is in the open arc $RQ_1$, and the Hessian is negative at the point $A\in L_1$ vertically above $A_2 \in L_2$.   Here we define $\bar A_2$ to be the point of the branch $B_2B_3$ of the Hessian immediately above $A_2$, 
i.e. lying on the line segment $A_2A$.  There is then a unique point $\bar U_2$ in the open arc 
$U_1 R$ of the Hessian with the tangent $\bar L_2$ at $\bar U_2$ containing  $\bar A_2$.  In passing, we remark  that previous arguments show that the region of $V_2$ given by $\bar A_2 \cdot D^2 >0$ consists only  of points below 
the line joining $\bar U_2$ and $\bar U_2' = \alpha (\bar U_2)$.  We denote by $Q$ the subcone of $P^\circ$ given by $\bar A_2 \cdot D^2 >0$; since $\tilde B_1\cdot D^2 >0$ for all $D \in P^\circ$, we note that $Q(2)$  is contained  in $Q$.  Since now both $U_1$ and $\bar U_2$ are contained in the open arc $Q_2R$, we have seen in 
the second of the above special cases that $Q(1)$ and $Q$ are disjoint, and hence $Q(1)\cap Q(2)$ is empty as required.
\end{proof}
\end{prop}

 \begin{cor}  Suppose that the real elliptic curve has $-2 < k <0$, then the statement 
 of Proposition 1.10 holds.
 
 \begin{proof} The argument here using limits may be reconstructed from the proof of Corollary 2.3, using the result 
 we have just proved, and is left as an exercise for the reader.
 \end{proof}
 \end{cor}
 
 \begin{figure}
 \centering
     \includegraphics[width=10cm]{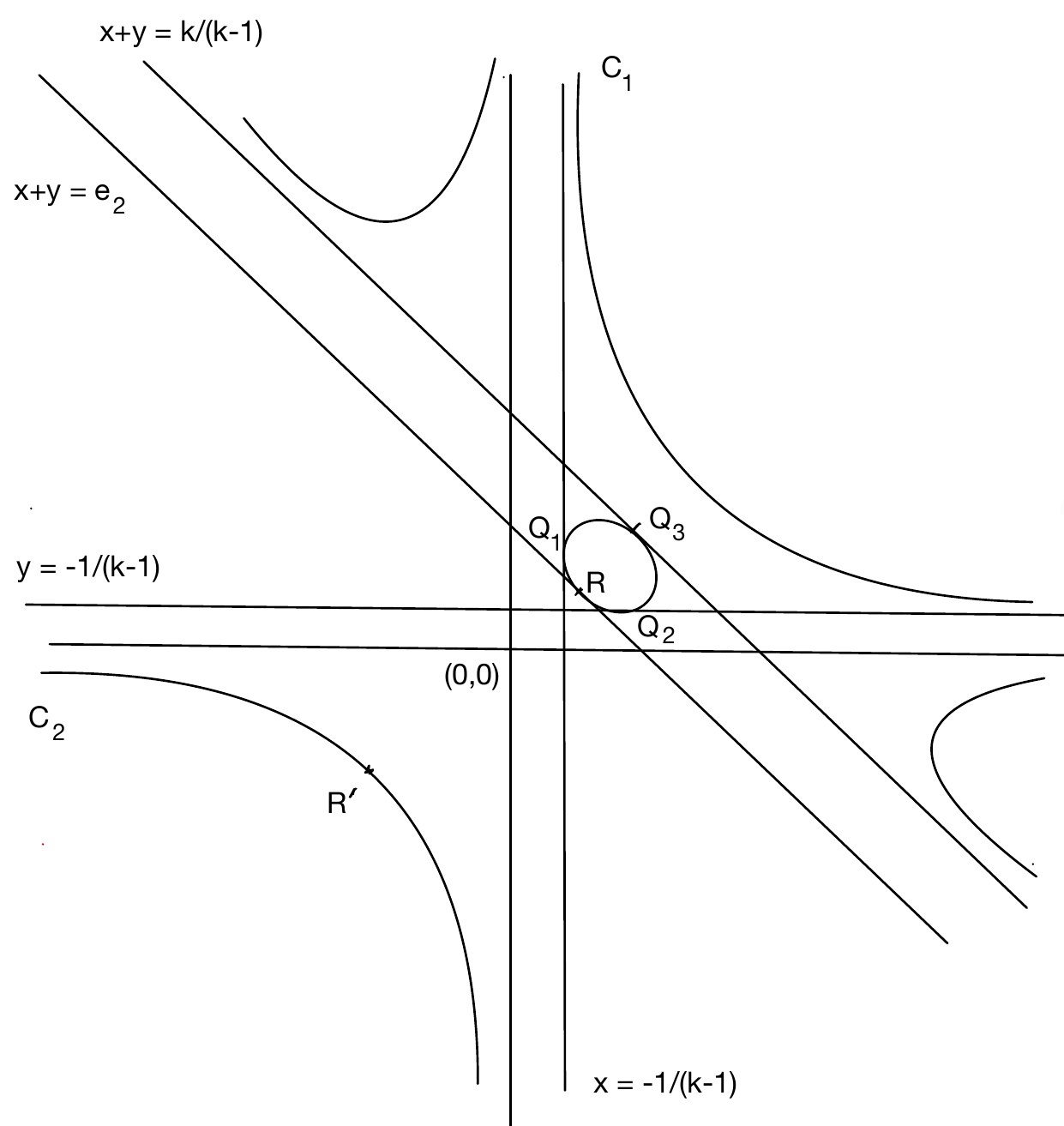}
\caption{Schematic picture when $k<-2$}
%\label{}
\end{figure}

\section{The case of the elliptic curve having invariant $k<-2$}

Let us now consider the remaining possibility with smooth Hessian, 
namely $k<-2$; 
here the bounded component of the elliptic curve is tangent to the asymptotic line $x+y = k/(k-1)$ but in this case lies below the line. We illustrate this situation schematically by Figure 6, where again we show the affine branches of the Hessian together with the affine branch $C_1$ of the cubic.  As usual, we assume that the hybrid component $P^\circ$ under consideration  has projectivised boundary 
corresponding to $C_1 \cup C_2$.   If the index at a point 
$E$ in the open upper half-space 
is $(1,q)$ with $q\le 2$, 
we  let $A=(a,b,1)$ denote the point in the affine plane $z=1$ determined by $E$.
  Let $Q$ denote a connected component of the subcone of $P^\circ$ given by $E\cdot D^2 >0$, 
by Proposition 3.4 of \cite{WilBd}, a convex subcone of $P^\circ$.  We now list the points on $C_1$ and $C_2$ where the quadratic $G_A$ vanishes.

As before it is clear for $A = (a,b,1)$ how many times (and where) $G_A=0$ intersects $C_1$.  If $a \ge {1\over {1-k}}$ 
and $b \ge {1\over {1-k}}$, then either $A\in V_1$ and so $G_A >0$ on both $C_1$ and $P^\circ$, or $G_A =0$ intersects $C_1$ twice (possibly at infinity or with multiplicty); it intersects $C_1$ 
 at no points if $a < {1\over {1-k}}$ 
and $b < {1\over {1-k}}$, and at precisely one point otherwise.  In the first of these three cases, we note that all of $C_2$ is visible from $A$.  
For 
$a < {1\over {1-k}}$, $b < {1\over {1-k}}$, there are no zeros of $G_A$ on  $C_2$ either, and a previous continuity argument shows that $G_A$ would be negative on all of $P$, since 
this is true for $A\in V_2$.  Thus if 
$a \le {1\over {1-k}}$, $b \le {1\over {1-k}}$, then $G_A$ would be negative on $P^\circ$ and so the subcone defined by $G_A >0$ is empty.  If 
for instance $a > {1\over {1-k}}$, $b = {1\over {1-k}}$ then $G_A$ has a zero on the affine branch $C_1$ and at  the point $B_2$ at infinity, and  
the whole of $C_2$ is visible from $A$.  We note that $G_A(B_1) >0$ if and only if $a > {1\over {1-k}}$ and $G_A(B_2) >0$ if and only if $b>{1\over {1-k}}$.

An additional feature compared with the previous case is that $G_A =0$ can intersect the projectivised boundary of $P$ at two points on $C_1$ and two on $C_2$, and that will happen when $A$ is in the open region with boundary consisting of a segment of the line $x = {1\over {1-k}}$, a segment of the line  
$y = {1\over {1-k}}$, and the `lower' arc (i.e. not containing $Q_3$) of the bounded component of the Hessian between $Q_1$ and $Q_2$; as before $Q_i$ denotes the point on the bounded component of the Hessian where the asymptote to the cubic at $B_i$ is tangent. In contrast to the previous case, this time the 
 arc 
of the Hessian between $Q_1$ and $Q_2$ lies in the quadrant $x \ge {1\over {1-k}}$, $y\ge {1\over {1-k}}$.  If $A$ lies in the quadrant  $x > {1\over {1-k}}$, $y> {1\over {1-k}}$ but not 
on or below the arc $Q_1Q_2$, then $G_A$ is non-zero (and therefore positive) on $C_2$.  
We note that the possibilities $E= \pm \tilde B_i$ for $i=1,2$ will give rise to trivial cases where either $G_E >0$ on both $C_2$ and $P^\circ$, or $G_E <0$ on $P^\circ$.
We remark in passing that if $-E$ is in the interior of the bounded component of the Hessian (and hence $E$ has $E^3 >0$ and index $(1,2)$, and is in the lower half-space), then $G_E$ is negative on $C_2$ but $P^\circ \cap \{ G_E >0 \}$ is non-empty; such an eventuality does not however occur under our assumptions.

We now need to understand what happens when $A$ lies on the tangent line at some point $U$ in the 
lower arc $Q_1Q_2$ of the bounded component of the Hessian.  We assume by symmetry that $U$ lies in the open arc $Q_1 R$, leaving it to the reader to check what happens when $U$ is $Q_1$ or $R$, including the limit case when $E$ is a positive multiple of $\pm \tilde B_3$.  
Under the given assumption, the tangent line meets the branch $B_1B_3$ of the Hessian.  We note that $A = (a,b,1)$ is below (and to the right of) 
this intersection point.  
If $a < {1\over {1-k}}$,  we know that $G_A$ vanishes at $U'$ and is positive at $B_2$, and there is 
a unique component of $P^\circ \cap \{ G_A > 0 \}$, whose closure contains $\tilde B_2$; the second point 
on the projectivised boundary of $P$ corresponds to the point on $C_1$ where the tangent line contains $A$.
Thus the  set of points in $C_2$ whose negative multiples are in the boundary
 of $P^\circ \cap \{ G_A > 0 \}$  is just the arc 
$B_2U'$, where as usual $U'=\alpha (U)$.
When $a = {1\over {1-k}}$, we know that $G_A$ vanishes (twice) at $B_1$, and for $A$ between this point and $U$, there is a point $U_1$ on the arc $Q_1U$ whose tangent also contains $A$, and there are two components of $P^\circ \cap \{ G_A > 0 \}$, one of which has its boundary points corresponding in $C_2$ to 
the arc $B_2 U'$ and the other with boundary points corresponding in $C_2$ to the arc $U_1'B_1$.  For 
$A= U$, we get two real lines meeting at $U'$, each line joining $U'$ to one of the two points of $C_1$ for which the tangent contains $U$.  
Moving now to $A$ lying below $U$ on the tangent  but with $b > {1\over {1-k}}$, we obtain a point $U_2$ in the arc $UQ_2$ 
where the tangent contains $A$, and we still have two components, but with the relevant arcs on $C_2$ 
now being $B_2 U_2'$ and $U' B_1$.  By the time we reach the point with $b = {1\over {1-k}}$, the first of these components 
of $P^\circ \cap \{ G_A > 0 \}$ has shrunk to the empty set (with $G_A$ vanishing along the ray generated by $\tilde B_2$), 
and from then on we just have one component of the intersection, and the relevant arc in $C_2$ is now 
$U'B_1$.  With this description in hand, the required result on visibility follows easily.

  \begin{prop}  Suppose the elliptic curve has one component, with invariant $ k<-2$.
  Given a component $Q$ of $P^\circ \cap \{ E\cdot D^2 >0 \}$, associated to a real  class $E$ in the  extended upper half-space  with 
 index $(1,q)$ for $q\le 2$,
  all points of $C_2$ whose negative multiples are on the boundary of $Q$ are visible 
  (with respect to the cone $-Q$) from $E$.

\begin{proof} 
Since visibility is a closed property, we may assume that $E$ is represented by a point $A$ in the affine plane.  
There were two general cases noted above where the result was true for trivial reasons --- in all other cases, $A$ lies on the tangent line at $U$ for some $U$ in the lower arc $Q_1Q_2$ on the bounded component of the Hessian.

We now just observe, in the various possibilities for $A$ lying on the tangent line at $U$, which points of $C_2$ are 
visible from $A$.  By symmetry (and noting that visibility is a closed condition) we may assume that $U$ lies in the open arc $Q_1R$.
 By the classical result used in previous sections, the point of intersection of the tangent at $U$ with the branch $B_1B_3$ of the Hessian is also just the intersection with the tangent line to $C_2$ at $U'$.  Thus the whole closed arc $B_2 U'$ is visible from $A$.  As $A$ moves down the tangent line, the 
visible points from $A$ constitute a larger arc containing the arc $B_2U'$, until  we reach the point 
with $a = -1/(k'-1)$  where $k' = (4-k^3)/(3k^2) >1$, when all of $C_2$ is visible from $A$.  This then remains true until $A$ has 
$b < -1/(k'-1)$, at which stage the visible points of $C_2$ nonetheless  form an arc containing $U'B_1$.  Thus we have 
verified the Proposition in all cases.
\end{proof}
\end{prop}

 \begin{prop}  Suppose the elliptic curve has one component, with invariant $ k<-2$.
 Let $E_i$ (for $i=1,2$) be real classes of UHS type with 
 indices $(1,q_i)$ for $q_i\le 2$,
  and suppose 
  there are  components $Q(i)$ of $P^\circ \cap \{G_{E_i}>0\}$ for $i=1,2$ whose intersection is non-empty; then some non-trivial  open arc in $C_2$ is in the boundaries of 
 both $-Q(1)$ and $-Q(2)$.

\begin{proof}
We show  that if the arcs in $C_2$ corresponding to the $Q(i)$ are disjoint, then $Q(1) \cap Q(2)$ is empty. We assume that neither $E_i$ correspond to a case where the result is trivially true. 
   Analogous to the proof of Proposition 3.2, 
 we then have points $U_1$, $U_2$ on the lower arc 
 $Q_1 Q_2$ of the bounded component (with $U_2$ in the arc $Q_1 U_1$), with corresponding 
 images $U_1'$, $U_2'$ on $C_2$ (with $U_2'$ in the arc $U_1'B_1$), where without of generality the open arc 
 $\Gamma _1 \subset C_2$ corresponding to $Q(1)$ is the arc $B_2U_1'$ and 
 the open arc $\Gamma_2 \subset C_2$ corresponding to $Q(2)$ is the arc $U_2'B_1$.
 As in the proof of Proposition 3.2, we shall assume that the $E_i$ define points $A_i$ in the affine plane, 
 since as we note below the limit cases $E_i = \pm \tilde B_3$ follow by only a 
 minimal change to the argument.  
 Thus $Q(1)$ corresponds to a point $A_1$ on the tangent line $L_1$ at $U_1$, with $A_1$ above (and to the left of) $U_1$, whilst $Q(2)$ corresponds to a point $A_2$ on the tangent line $L_2$ at $U_2$, with $A_2$ below (and to the right of)  $U_2$. Thus $Q(1)$ contains $\tilde B_2$ and 
 $Q(2)$ contains $\tilde B_1$ in their boundaries.  
 Under these conventions, we show that $Q(1)\cap Q(2)$ is empty.  
 
 In the above notation, $L$ is just the line joining $U_1$ and $U_1'$, whilst $M$ is the line joining $U_2$ and $U_2'$.  We observe that as $U$ moves monotonically from $Q_1$ to $Q_2$ on the arc $Q_1Q_2$ of the Hessian, its conjugate $U' = \alpha (U)$ moves monotonically from $B_1$ to $B_2$ on $C_2$.  We deduce that $L$ and $M$ intersect at an affine point not in $V_1 \cup V_2$, and so the corresponding planes through the origin 
 in $\R^3$ intersect in a line disjoint from $P^\circ$.  The plane determined by $L$ gives rise to a cone $\tilde Q(1)$, which is the intersection of the relevant open half-space containing $\tilde B_2$ with $P^\circ$, and the plane determined by $M$ gives rise to a cone $\tilde Q(2)$, which is the intersection 
  with $P^\circ$  of the relevant  open half-space containing $\tilde B_1$.  Thus $Q(1) \subset \tilde Q(1)$ and $Q(2) \subset \tilde Q(2)$, which by the usual continuity argument from Section 2 remains true even if one or both of the $A_i = U_i$.
 Since $\tilde Q(1) \cap\tilde Q(2)$ is empty, we obtain the claimed result.
 \end{proof}
\end{prop}

\begin{cor}  Suppose that the real elliptic curve has $k < -2$, then the statement 
 of Proposition 1.10 holds.
 
 \begin{proof} The argument via limits may again be reconstructed from the proof of Corollary 2.3, using the result 
 we have just proved, and is left as an exercise for the reader.
 \end{proof}
 \end{cor}
 
 Therefore in the case when the real elliptic curve has one component, we have completed (via 
 Proposition 1.9 and Corollaries 3.3 and 4.3) the proof of Theorem 1.8,  and hence 
we have proved (via Proposition 1.4) the relevant part of our
Theorem 0.2.

\section{Boundedness questions}

Finally we prove Theorem 0.3.

\begin{proof}  The assumption that $E^3 \le 0$ for all rigid non-movable surfaces $E$ on $X$ implies that $c_2(X)\cdot E \ge 0$ for all such surfaces by Proposition 2.2(i) of \cite{WilBd}.  
We saw in Corollary 1.7 that 
the integral linear form $c_2$ is non-negative on the closed positive cone on the bounded component of the real elliptic curve.  Having chosen \it real \rm
coordinates as in equation (1), we let $S=S^2$ denote the unit sphere in $\R^3$; intersecting the plane $c_2 =0$ with the upper half-sphere, we obtain a half-circle with end-points on the equator.  
Our assumptions imply that the cubic is strictly positive at one of these endpoints $B$, and then 
there exists a point $T$ on the half-circle with the cubic strictly 
positive on the whole of the  arc $TB$. 
   If the line $c_2=0$ is tangent to the bounded component of the cubic, then note that $T$ lies between $B$ and the point on the half-circle corresponding to  the point of tangency.

Since $c_2$ is non-negative on the closed bounded component of the elliptic curve (topologically a closed disc), the map $L\mapsto L^2$ identifies the boundary of the bounded component with the boundary of its polar (with a boundary point corresponding to its tangent line), with the polar also topologically a closed disc.  
If the map is not surjective, then there is 
a continuous retraction of the disc onto $S^1$, 
a contradiction;  the corresponding map on the positive cone on the bounded component therefore identifies the dual cone as the image of this map,
 The map is however also injective since if $D_1 ^2 = D_2^2$, then $(D_1+D_2)(D_1 -D_2) =0$ and 
the Hessian at $D_1+D_2$ is non-zero by convexity, implying that $D_1 = D_2$.
There therefore exists a unique ray $\R _+ D_0$ with $D_0$ an affine real point in the bounded component with the linear form $c_2$ being a positive multiple of $D_0^2$; the case when $D_0$ is in the boundary  corresponds to the line $c_2 =0$ being tangent to the bounded component.
If $A_1, A_2$ are points on the projective  line $c_2 =0$, then $D_0$ is the point in the bounded component given by intersecting the conics $A_i \cdot D^2 =0$ for $i=1,2$.  Now choose rational points $A_1, A_2$ on the projective line $c_2 =0$ 
corresponding to points in  the open arc $TB$, therefore not lying on the Hessian  
and hence  defining  smooth projective conics $G_{A_1}$ and $G_{A_2}$.  By taking an appropriate rational convex combination $A$ of $A_1$ and $A_2$, the conic $G_A$ is the corresponding rational combination of $G_{A_1}$ and $G_{A_2}$, which we may arrange to vanish at a rational point; since $G_A$ is also smooth, its rational points are then dense.

By taking a rational point $D$ on this conic sufficiently close to $D_0$ and in the interior of the bounded component, we may arrange that 
$D^2 =0$ is a projective line through $A$, defining a 
half-circle in $S$ which intersects the equator at a point $B'$ near $B$, with the linear form $D^2$ being negative at $B$, for which the cubic is strictly positive on  the non-zero points of the closed simplicial cone generated by $A$, $B$ and $B'$.  We deduce that there are no 
potential rigid non-movable classes in this closed cone, since by assumption any such class $E$ has $E^3\le 0$.  
Since the class of any rigid non-movable surface lies in the extended upper half-space, it 
follows then that $D^2 \cdot E >  0$ for all rigid non-movable classes $E$ on $X$.  Without loss of generality, we take the class $D$ to be integral and 
assume  that $X$ is general in moduli; we now apply Proposition 4.1 and Lemma 4.3 from \cite{WilBd} to produce an integer $m>0$, depending on $D$ but not on the rigid non-movable classes $E$, such that $h^0 (X, \O_X (nD)) >1$ for all $n\ge m$, with moreover $D$ of general type.  We can then write $|mD| = |\Delta | + \calE$, with $\calE$ a finite sum $\Sigma  r_i E_i$ supported on rigid non-movable surfaces $E_i$ (with each $r_i>0$) and $\Delta$ mobile.

For any mobile class $M$, we know that $c_2\cdot M \ge 0$; we now claim also that $D^2 \cdot M \ge 0$.  To see this, suppose first that $M$ lies in the closed lower half-space --- the fact that the Hessian is non-negative at $M$ then implies that $M^3 \ge 0$.  Choose a very ample divisor $L$; the divisors $L + tM$ are then mobile with strictly positive cube for all $t\ge 0$, which implies that $M$ is in the closure of the component of the positive index cone containing $L$, namely the positive cone on the bounded component of the elliptic curve, and this is clearly nonsense (cf. the last paragraph in the proof of Proposition 1.9, the last part of the proof of Proposition 1.11 and the argument in Remark 1.12).  Thus 
a mobile class $M$ must lie in the open upper half-space.  If however it lies in the closed simplicial cone generated by $A$, $B$ and $B'$, we know that 
$M^3 >0$ and a similar contradiction is obtained.  Therefore the fact that $c_2\cdot M \ge 0$ now implies that $D^2 \cdot M \ge 0$..

By considering $mD\cdot c_2$ and $mD^3$, we note that for any $i$
we have bounds on $c_2 \cdot E_i \ge 0$, $D^2\cdot E_i >0$ and $r_i$.   With respect to the standard norm on $\R^3$,  
we claim that $\|E \|$ is bounded for all possible classes $E$ which might occur as an $E_i$ in a decomposition of $mD$, and hence by compactness there are only 
finitely many possibilities for these classes.  
If not, then there are infinitely many classes $F_j$ which could 
potentially represent such an $E_i$, therefore all lying in the extended upper half-space and with $c_2 \cdot F_j$ and $D^2\cdot F_j$  bounded, but with $\| F_j\|$ unbounded.  
Thus some subsequence of the points $\bar F_j := F_j /\| F_j \|$ on the sphere $S$ has a limit point in $S$, which must lie in the upper half-sphere.
Since both $c_2 \cdot \bar F_j $ and $D^2 \cdot \bar F_j $ tend to zero for this subsequence, this limit point can only be the point $A/\|A\|$;
this is a contradiction since $A$ was chosen so that $A^3 >0$, whilst by assumption  $E_i^3 \le 0$ for all $i$.  As there are only finitely many possibilities for the classes $E_i$, 
and the $r_i$ are bounded, 
there are only finitely many possible classes for the class of the mobile part $\Delta$ of $|mD|$. We may assume therefore that the class of the mobile divisor $\Delta$ is also given; 
 under the assumption that the divisor is big, we may deduce boundedness by Proposition 1.1 of \cite{WilBd}.  

Otherwise, 
it may be that  the mobile divisor $\Delta$ defines a fibration stucture on some minimal model of $X$.  The case when the fibration structure is a fibre space over $\P^1$ is easy 
to deal with, since if we write $\Delta = NG$ for some (known) $N$ and primitive class $G$, then on the appropriate minimal model $\tilde X$ the fibration corresponds to the semi-ample  
irreducible divisor $\tilde G$ corresponding to $G$, 
and it follows that $h^0 (X, \O_X (r\Delta)) = rN+1$ for all positive integers $r$.  Note that the divisor $\calE$ cannot contain the primitive class $G$, that is we cannot have $G\le \Delta$ as divisor classes.
Thus $h^0 (X, \O_X (m\Delta)) = mN+1$;  moreover for $n\ge 1$,  $h^0 (X, \O_X (nmD))$
is bounded below by $\frac{1}{6} n^3m^3 D^3 + \frac{1}{12} nm c_2\cdot D$ by Proposition 4.1 and Lemma 4.3 from \cite{WilBd}.  For such an $n$,  if $n\calE$ contains some multiple $sG$ of the primitive class $G$, then $s <n$.  Thus if the mobile part of $|nmD|$ were to correspond to a fibration over $\P^1$ on the minimal model, then $h^0 (X, \O_X (nm\Delta)) \le nmN+n$.  We can however choose $n$ so that $
\frac{1}{6} n^3m^3 D^3 + \frac{1}{12} nm c_2\cdot D > nmN + n$, and  we deduce that 
the mobile part of $| nmD |$ does not correspond to a fibration over $\P^1$ on the minimal model.  

We are therefore left with the cases when some (known)  mobile divisor  on $X$  defines an 
elliptic fibre space structure on some minimal model 
and so $X$ is birationally elliptic.  In particular, using the results from \cite{Gr}, we know  that the given family of \CY threefolds is at least birationally bounded.
\end{proof}

\begin{rem}  
Under the assumptions of the theorem,  we deduce that the third betti number $b_3$ is bounded by knowledge of the cubic and linear forms on 
the integral second cohomology.  Thus for a given smooth integral 
ternary cubic with irreducible Hessian, and a linear integral form $c_2$ not vanishing at any real inflexion point of the elliptic curve,  
there exists an integer $b\ge 0$ (depending on the forms) such that any \CY threefold $X$ with  $b_3 (X) > b$ which realises the given cubic and linear forms  must contain  a rigid non-movable surface $E$ with $E^3 >0$ (and by Remark 1.12 at most six of them); in particular, after at most eight flops on $X$ we obtain a minimal model on which the transform of $E$ may be contracted
\end{rem}

% ----------------------------------------------------------------
\bibliographystyle{amsalpha}

\end{document}